\documentclass{article}[12pt]
\usepackage{amssymb}
\usepackage{epsfig}
\usepackage{amsmath}
\usepackage{amsfonts}
\usepackage{amsthm}
\usepackage{pgf,tikz}
\usepackage{enumerate}
\usepackage{cite}
\usepackage[colorlinks=true, linkcolor=black, citecolor=black, urlcolor=black]{hyperref}

 \usepackage{soul} 

\newtheorem{theorem}{\sc Theorem}[section]
\newtheorem{proposition}[theorem]{\sc Proposition}
\newtheorem{lemma}[theorem]{\sc Lemma}
\newtheorem{definition}[theorem]{\sc Definition}
\newtheorem{remark}[theorem]{\sc Remark}
\newtheorem{corollary}[theorem]{\sc Corollary}
\newtheorem{example}[theorem]{\sc Example}

\def\qed{\hbox to 0pt{}\hfill$\rlap{$\sqcap$}\sqcup$\medbreak}

\title{Lower and upper solutions for boundary value problems with bounded $\Phi$-laplacian Stieltjes-differential inclusions}

\author{Rodrigo L\'opez Pouso and Paula Vázquez García}

\date{}
\begin{document}
 \maketitle

\begin{center}  {CITMAga \& Departamento de Estat\'{\i}stica, An\'alise Matem\'atica e Optimizaci\'on, \\ Universidade de Santiago de Compostela, \\ 15782, Facultade de Matem\'aticas, Campus Vida, Santiago, Spain.}
\end{center}

\medbreak

 \abstract{We prove the existence of solutions between given lower and upper solutions for several types of boundary value problems involving bounded $\Phi$-Laplacian inclusions. We consider derivatives in the Stieltjes sense, that is, derivatives with respect to nondecreasing and left-continuous functions. This framework includes, as particular cases, both differential and difference equations, even allowing for impulses at the discontinuity points of the derivator. To the best of the authors’ knowledge, this paper presents the first application of the method of lower and upper solutions to second-order boundary value problems with Stieltjes derivatives.}

\medbreak


\medbreak


\section{Introduction}\label{sec1}
According to Esteban and Vázquez \cite{esva}, we owe to Leibenson \cite{leib} the following one--dimensional model for the turbulent flow of a polytropic gas in a porous media:
\begin{equation}
\label{pde}
\dfrac{\partial u}{\partial t}= k \, \dfrac{\partial}{\partial x} \left( 
\dfrac{\partial u^m}{\partial x} \left| 
\dfrac{\partial u^m}{\partial x} \right|^{p-1} \right).
\end{equation}
Here, $u$ is the density of the gas and the constants $k$, $m$ and $p$ are positive. This type of nonlinear diffusion equation has drawn a lot of attention in the 1980s, and we can quote many relevant papers such as \cite{alieva, esva, heva} and references therein. If stationary solutions of (\ref{pde}) are considered and the assumption is made that the gas is being consumed by an irreversible reaction, then we end up with ordinary differential equations of the form
\begin{equation}
\label{ode}
\dfrac{d}{d x} \left( 
\dfrac{d u^m}{d x} \left| 
\dfrac{d u^m}{d x} \right|^{p-1} \right)=f(x,u,u'),
\end{equation}
where the choice of $f(x,u,u')$ depends on the reaction. When $m=1$, the differential operator in the left--hand side of (\ref{ode}) is called the $p$--laplacian (which reduces to the usual laplacian $u''$ if $p=1$). Bobisud \cite{bob} observed that (\ref{ode}) could be studied as a particular case of an equation
\begin{equation}
\label{fibo}
(\Phi(u'))'=f(x,u,u'),
\end{equation}
where $\Phi: \mathbb R \longrightarrow \mathbb R$ is a given homeomorphism, starting in this way the study of the so--called $\Phi$--laplacian equations. See \cite{capo, capo2, capo3, capo4, cahapo, ore, wgl} and references therein. For consistency with the notation in these and subsequent references (and also with the notation in this paper), we replace the $x$ variable with the $t$ variable and rewrite (\ref{fibo}) as
\begin{equation}
\label{fibot}
(\Phi(x'))'=f(t,x,x'),
\end{equation}
where now $x(t)$ is the unknown solution.

Bereanu, Jebelean and Mawhin studied in a series of papers \cite{bema, bema2, bema3, bema4, bejema} the case of possibly bounded $\Phi$ in (\ref{fibot}) and the case of a function $\Phi$ defined on a bounded interval. Model examples are
$$\Phi(z)=\dfrac{z}{\sqrt{1+z^2}}, \, \, z \in \mathbb R, \quad (bounded)$$
and 
$$\Phi(z)=\dfrac{z}{\sqrt{1-z^2}}, \, \, z \in (-1,1), \quad (singular)$$
which correspond to the mean curvature operator in Euclidean space and in Minkowski space, respectively; see \cite{bejema}.

Yet another motivation for the study of $\Phi$--laplacian equations of the form of (\ref{fibot}) is indicated in \cite{capo4} in connection with right--hand sides $f(t,x,x')$ which are not continuous. Due to insufficiently developed theory, discontinuities in differential equations are often replaced with all sorts of regularizations, but this is not satisfactory when the presence of discontinuities is natural. Consider, for instance, a simple model of the damped motion of a stone falling from the air into the water, which clearly demands a discontinuous damping term to reflect the sudden change of medium. As we pointed out in \cite{capo4}, any ordinary differential equation of the form
\begin{equation}
    \label{capoedo}
x''=f(t,x,x')q(x'),
\end{equation}
where $f(t,x,x')$ is continuous and $q, \, 1/q \in L^{\infty}_{loc}(\mathbb R,(0,\infty))$, can be equivalently expressed as (\ref{fibot}) with $$\Phi(z)=\int_0^z{\dfrac{dv}{q(v)}},\quad z\in\mathbb{R}.$$
Hence, the discontinuous ODE (\ref{capoedo}) is transformed into an equivalent continuous $\Phi$--laplacian equation.

Similarly, we can motivate the introduction of Stieltjes derivatives in (\ref{fibot}). The Stieltjes derivative of a real function $F(t)$, see \cite{FeMaTo, FrLo, LoMa, LoRo, mosa}, with respect to a nondecreasing and left--continuous function $g: \mathbb R \longrightarrow \mathbb R$ (which we call {\it derivator}), is essentially defined as
\begin{equation}\label{derg}
    F_g'(t_0)= \lim_{t \to t_0^+}\dfrac{F(t)-F(t_0)}{g(t)-g(t_0)},
\end{equation}
unless $g$ is constant on $[t_0,t_0+\varepsilon)$ for some $\varepsilon>0$. Stieltjes differential equations include as particular cases both ordinary differential equations and difference equations, even with impulses at fixed times; see \cite{FrLo, LoMa, LoRo}. Here we exploit an interesting result by Monteiro and Satco \cite{mosa}, which ensures that, under some technical but very general conditions, a distributional differential equation
$$D(x')=f(t,x,x')Dg,$$
where $D(x')$ and $Dg$ stand for the distributional derivative of $x'$ and $g$, respectively, can be seen as a particular case of a Stieltjes differential equation
$$(x')_g'=f(t,x,x'),$$
where $(x')'_g$ is the (pointwise) Stieltjes derivative of $x'$. In doing so, we pass from a complicated discontinuous distributional equation to a simpler pointwise-defined Stieltjes differential equation. To have a clearer idea of how useful these reductions are, consider the singular Thomas--Fermi equation
\begin{equation}
    \label{stf}
x''=x^{3/2}t^{-1/2}, \quad t>0,
\end{equation}
which can be rewritten as
\begin{equation}
    \label{stfg}
(x')_g'=x^{3/2}, \quad t>0,
\end{equation}
where $g(t)=2 \sqrt{t}$ for $t  \ge 0.$ Observe that equation (\ref{stf}) is singular at $t=0$ but equation (\ref{stfg}) is not.

\medbreak

This paper is devoted to the study of different types of boundary value problems for differential inclusions of the form
\begin{equation}\label{eqg}
   ((A\circ x) \cdot(\Phi \circ x'_{g_1}))'_{g_2}(t) \in G(t,x(t),x'_{g_1}(t)), \quad t \in I=[a,b],
\end{equation}
where, besides some other technical conditions to be detailed in the next sections, we assume that
\begin{itemize}
\item $g_1, \, g_2: \mathbb R \longrightarrow \mathbb R$ are nondecreasing and left-continuous everywhere;
    \item $A:\mathbb{R}\longrightarrow\mathbb{R}$ is a positive and continuous function;
    \item $\Phi:\mathbb{R}\longrightarrow(-c,c)$, $c>0$, is an increasing homeomorphism; and
    \item $G:I\times\mathbb{R}\rightarrow\mathcal{P}(\mathbb{R})$ is a Carath\'eodory multivalued mapping, with non-empty, compact and convex values.
    \end{itemize}

We follow the methodology employed by Ferracuti, Marcelli and Papalini \cite{FeMaPa}, where the same type of differential inclusions is considered using standard derivatives. The introduction of Stieltjes derivatives gives rise to several technical complications, which the reader can readily appreciate by comparing our proofs in Section 4 with those in \cite{FeMaPa}.

This paper is organized as follows. In Section \ref{sec2}, we include all the necessary technical preliminary results on set-valued analysis and Stieltjes calculus; we also define in Section \ref{sec2} two concepts of solutions of a differential inclusion, namely weak and strong solutions. In Section \ref{sec3}, we prove the existence of weak and strong solutions for a technical problem by means of which we can study the general problems \eqref{eqg}. In Section \ref{sec4}, we prove our main existence result for \ref{eqg} with Dirichlet boundary conditions and lower and upper solutions. To the best of the authors' knowledge, this is the first application of the method of lower and upper solution with second order equations of Stieltjes type. In Section \ref{sec4} we also prove a lemma on continuous dependence of the solutions with respect to the boundary conditions. This lemma is crucial in the proof of the results in Section \ref{sec5}, where we consider nonlinear boundary conditions which include periodic and Robin boundary conditions as particular cases. Examples of applications of our theorems are provided.

\section{Preliminaries}\label{sec2}
\subsection{Set-valued analysis}
In this section we gather some important definitions and results on multivalued mappings. In the sequel, $X$ and $Y$ denote two normed spaces.

\begin{definition}
    A multivalued mapping $F:X\to\mathcal{P}(Y)$ is
    \begin{enumerate}
        \item \textit{upper semicontinuous} (u.s.c.) at a point $x\in X$ if for each open neighborhood $U$ of $F(x)$ there exists $\varepsilon>0$ such that 
    $F(y)\subset U$ for every $y\in B_\varepsilon(x)$. We say that $F$ is upper semicontinuous if the condition is satisfied at every point of its domain.
        \item \textit{weakly upper semicontinuous} (w-u.s.c.) at a point $x\in X$ if for each sequence $\left\{x_n\right\}\subset X$ such that $x_n\to x$ and each sequence $\left\{y_n\right\}$ such that $y_n\in F(x_n)$ and $y_n\to y$ weakly in $Y$ we have that $y\in F(x)$. 
    \end{enumerate}
\end{definition}

The graph of a multivalued mapping $F:X\to\mathcal{P}(Y)$ is defined as the set
$$
    \mathrm{Gr}(F) = \{ (x,y) \in X \times Y : y \in F(x) \}.
$$
If the graph is closed, and $F$ is sufficiently nice, then $F$ is u.s.c. More precisely, we have the following.

\begin{proposition}
    \label{prohupa}
    \cite[Proposition 2.23]{HuPa1}
    Suppose that $F:X\to\mathcal{P}(Y)$ assumes nonempty closed values and is locally compact (i.e., every $x \in X$ has a neighborhood $U$ such that $\overline{F(U)}$ is compact).
    
      If $F$ has a closed graph, then $F$ is u.s.c.
\end{proposition}

\bigbreak

We will use the following general fixed point theorem proven by Ferracuti and Papalini \cite{fepa} as a generalization of an earlier result by Cellina \cite{cel}. Remarkably, it ensures the existence of at least one fixed point for not necessarily convex valued mappings.

\begin{theorem}\label{teo:fp1}
    Let $S$ be a non-empty, compact, and convex subset of a Banach space $X$ and $D$ be a non-empty, convex subset of a normed space $Y$. Let $F:S\to\mathcal{P}(D)$ be an u.s.c. multivalued mapping with non-empty and convex values such that $\mathrm{co}F(S)\subset D$, and $h:S\times D\rightarrow S$ be a uniformly continuous singlevalued function. Then, the multivalued mapping $\Gamma:S\rightarrow\mathcal{P}(S)$ given by $\Gamma(x):=h(x,F(x))$, $x\in S$, has at least one fixed point provided that $\Gamma$ has a closed graph.
\end{theorem}
\subsection{Background on Stieltjes calculus}
In this section, we gather some concepts and results on Stieltjes calculus which will be needed in the sequel. The Stieltjes calculus leans on the notion of a derivator: we say that $g:\mathbb R \longrightarrow \mathbb R$ is a derivator if $g$ is monotone nondecreasing and continuous from the left everywhere. We need the set of points around which $g$ is constant, that is,
\begin{equation}
\label{Cg}
C_g=\{ t \in \mathbb R \, : \, \mbox{$g$ is constant on $(t-\varepsilon,t+\varepsilon)$ for some $\varepsilon>0$}\},
\end{equation}
and the set of discontinuity points of $g$, which is at most countable and shall be denoted by
\begin{equation}
\label{Dg}
D_g=\{ t \in \mathbb R \, : \, g(t^+)-g(t)>0 \},
\end{equation}
where, as usual, $g(t^+)=\lim\limits_{s\to t^+}g(s)\in\mathbb R$.

Following \cite{LoRo}, we define the {\em derivative with respect to $g$} (or {\em$g$--derivative}) of a real--valued real function $f$ at a point $t_0 \in \mathbb R \setminus C_g$ as follows, provided that the corresponding limit exists:
\begin{equation}
\label{dercont}
f'_{g}(t_0)=\lim_{t \to t_0}\dfrac{f(t)-f(t_0)}{g(t)-g(t_0)} \quad \mbox{if $t_0 \not \in D_g$, or}
\end{equation}
\begin{equation}
\label{derdis}
f'_g(t_0)=\lim_{t \to t_0^+}\dfrac{f(t)-f(t_0)}{g(t)-g(t_0)} \quad \mbox{if $t_0 \in D_g$,}
\end{equation}
where $C_g$ and $D_g$ are defined in~\eqref{Cg} and~\eqref{Dg} respectively. Notice that the $g$--derivative at a point $t_0 \in D_g$ exists if, and only if, $f(t_0^+)$ exists and~\eqref{derdis} can be rewritten as
\begin{equation}
\label{derdis2}
f'_g(t_0)=\dfrac{f(t_0^+)-f(t_0)}{g(t_0^+)-g(t_0)}.
\end{equation}
We say that $f$ is {\em $g$--differentiable at $t_0$} if $f'_g(t_0)$ exists, and we say that $f$ is {\em $g$--differentiable in a set $A \subset \mathbb R$} when $f$ is $g$--differentiable at every $t_0
\in A \setminus C_g$. 

There is often no need to define $g$--derivatives for the points in $C_g$ because, according to~\cite[Proposition~2.5]{LoRo}, we have $\mu_g(C_g)=0$, where $\mu_g$ stands for the Lebesgue--Stieltjes measure induced by $g$, which satisfies $\mu_g([c,d))=g(d)-g(d)$ for every $c,d \in \mathbb R$, $c<d$. While Lebesgue-Stieltjes integration is analogous in most aspects to Lebesgue integration, there are some important differences. To start with, the Lebesgue-Stieltjes measure can have atoms, i.e., singleton sets with positive measure (which correspond to points in $D_g$). Specifically, we have $\mu_g(\{t\})=g(t^+)-g(t)>0$ when $t \in D_g$ (see equation~(6) in~\cite{LoRo}). 

We say that a property is satisfied for $g$-almost all $t\in A \subset \mathbb R$ or $g$-almost everywhere in $A$ ($g$-a.a. $t \in A$ or $g$-a.e. in $A$, for short) if it is satisfied for all $t \in A \setminus B$ for some $B \subset A$ such that $\mu_g(B)=0$.  

We denote by $L_g^1(J)$ the set of (equivalence classes of $g$-almost everywhere identical) Lebesgue--Stieltjes integrable functions in the interval $J$, and for $f \in L_{g}^1(J)$, we denote its Lebesgue--Stieltjes integral in $J$ by
$$\int_{J}f(s) d\mu_g \quad \mbox{or} \quad \int_{J}f(s) dg(s),$$
and $L^1_g(J)$ is a Banach space with the norm
$$\|f\|_{L^1_g}=\int_{J}{|f(s)|\, dg(s)}, \quad f \in L^1_g(J).$$
 
Roughly speaking, $g$--differentiable functions can be recovered by integrating their $g$--derivatives in the Stieltjes sense with respect to $g$. More precisely, this is true for $g$-absolutely continuous functions.

\begin{definition}
\label{defacg}
    Let $a,b \in \mathbb R$, $a<b$, and $F:[a,b] \longrightarrow \mathbb R$. We say that $F$ is {\em absolutely continuous with respect to $g$\/} (or {\em $g$--absolutely continuous\/}) if to each $\varepsilon>0$ there is some $\delta>0$ such that, for any family $\{(a_n,b_n)\}_{n=1}^{m}$ of pairwise disjoint open subintervals of $[a,b]$, the inequality
$$
\sum_{n=1}^{m}(g(b_n)-g(a_n))<\delta
$$
implies
$$
\sum_{n=1}^m|F(b_n)-F(a_n)|<\varepsilon.
$$
We denote by $\mathcal{AC}_g([a,b])$ the set of $g$-absolutely continuous functions in $[a,b]$.
\end{definition}

Indefinite Lebesgue-Stieltjes integrals are $g$-absolutely continuous.
\begin{proposition}\cite[Proposition 5.2]{LoRo}
\label{proindeflsg}
If $f \in L^1_g([a,b))$, then the function
$$F(t)=\int_{[a,t)}{f(s) \, dg(s)}, \quad t \in [a,b],$$
belongs to $\mathcal{AC}_g([a,b])$.
\end{proposition}

The converse of Proposition \ref{proindeflsg} is included in the Fundamental Theorem of Calculus for the Lebesgue--Stieltjes integral~\cite[Theorem~5.4]{LoRo}.

\begin{theorem}{\bf (Fundamental Theorem of Calculus for the Lebesgue--Stieltjes integral)}
\label{tf}
Let $a,b \in \mathbb R$, $a<b$, and $F:[a,b] \longrightarrow \mathbb R$. The following conditions are equivalent.
\begin{enumerate}
\item[(1)] $F \in \mathcal{AC}_g([a,b])$.
\item[(2)] The function $F$ fulfills the following properties:
\begin{enumerate}
\item[(a)] there exists $F'_g(t)$ for $g$--almost all $t\in [a,b)$ (i.e., for all $t$ except on a set of $\mu_g$ measure zero);
\item[(b)] $F'_g \in L^1_g([a,b))$, the set of Lebesgue--Stieltjes integrable functions with respect to $\mu_g$;
\item[(c)] for each $t  \in [a,b]$, we have
\begin{equation}
\label{barrow}
F(t)=F(a)+\int_{[a,t)} F'_g(s) \, d\mu_g.
\end{equation}
\end{enumerate}
\end{enumerate}
\end{theorem}

In particular, if $t \in [a,b) \cap D_g$, then,~\eqref{barrow} implies that
$$
F(t^+)-F(t)=\int_{\{t\}}F'_g(s) \, d\mu_g=F'_g(t) \mu_g(\{t\})=F'_g(t)(g(t^+)-g(t)),
$$
which may be nonzero. Therefore, $g$--absolutely continuous functions need not be continuous at discontinuity points of $g$. However, $g$--absolutely continuous functions have some nice properties. 

\begin{proposition}\cite[Proposition 5.3]{LoRo}
\label{bv}
If $F \in \mathcal{AC}_g([a,b])$, then it has bounded variation and is continuous from the left at every $x \in [a,b)$.

Moreover, $F$ is continuous in $[a,b] \setminus D_g$, where $D_g$ is the set of discontinuity points of $g$, and if $g$ is constant on some $(\alpha,\beta) \subset [a,b]$, then $F$ is constant on $(\alpha,\beta)$ as well.
\end{proposition}

For later convenience, we also include the following elementary observation.  

\begin{proposition}
\label{maxminac}
If $F, G \in \mathcal{AC}_g([a,b])$, then 
\begin{enumerate}
\item[(i)] $c F+d G \in \mathcal{AC}_g([a,b])$ for every $c,d \in \mathbb R$.
\item[(ii)] $|F| \in \mathcal{AC}_{g}([a,b])$.
\item[(iii)] If for each $t \in [a,b]$ we define
$$m(t)=\min\{F(t),G(t)\} \quad \mbox{and} \quad M(t)=\max\{F(t),G(t)\},$$
then $m, M \in \mathcal{AC}_g([a,b])$.
\end{enumerate}
\end{proposition}

\noindent
{\bf Proof.} The proofs of $(i)$ and $(ii)$ are straightforward from the definition of $g$-absolutely continuous function in the first part of Theorem \ref{tf}. Now, $(iii)$ is a consequence of $(i)$, $(ii)$, and the expressions
$$m(t)=\dfrac{F(t)+G(t)-|F(t)-G(t)|}{2} \quad \mbox{and} \quad M(t)=\dfrac{F(t)+G(t)+|F(t)-G(t)|}{2}.$$
\qed

The natural norm in $\mathcal{AC}_g([a,b])$ is
$$\|F\|_{\mathcal{AC}_g}=|F(a)|+\int_{[a,b)}{|F'_g(t)| \, dg(t)}=|F(a)|+\|F_g'\|_{L^1_g}, \quad F \in \mathcal{AC}_g([a,b]),$$
and $(\mathcal{AC}_g([a,b]),\| \cdot \|_{\mathcal{AC}_g})$ is a Banach space (observe that there is a linear isometry with the Banach space $\mathbb R \times L^1_g([a,b))$).

\bigbreak

This paper is concerned with second order ODEs, so we need second order derivatives. To have a second order derivative we first need an everywhere defined first order derivative and therefore the definitions \eqref{dercont} and \eqref{derdis} are not enough. To have everywhere defined first order derivatives, we follow Fernández, Márquez-Albés and Tojo \cite{FeMaTo}. First, we introduce some notation. Since $C_g$ is open, we have
\begin{equation}\label{cg}
C_{g}=\bigcup_{n\in\mathcal{C}}(a_n,b_n), 
\end{equation}
where $\mathcal{C}\subset\mathbb{N}$ and $(a_i,b_i)\cap(a_j,b_j)=\emptyset$ for every $i,j\in\mathcal{C}$ with $i\neq j$. In these conditions, we denote $N_{g}^{-}=\{a_n\}_{n\in\mathcal{C}}\setminus D_{g}$, $N_{g}^{+}=\{b_n\}_{n\in\mathcal{C}}\setminus D_{g}$ and $N_{g}=N_{g}^{-}\cup N_{g}^{+}$.
We are now in a position to introduce the definition of the $g$-derivative of a function defined at every point of its domain, first given in \cite{FeMaTo}.  
\begin{definition}
\label{defderevery}
     Let $g:\mathbb{R}\longrightarrow\mathbb{R}$ be a derivator such that $a\notin N_{g}^{-}$ and $b\notin C_{g}\cup N_{g}^{+}\cup D_{g}$. We define de the $g$-derivative of a function $f:[a,b]\subset\mathbb{R}\longrightarrow\mathbb{R}$ at a point $t\in [a,b]$ as
     \begin{equation}\label{fg}
         f'_{g}(t)=
         \left\lbrace
         \begin{array}{cll}
            \vspace{3mm}
             \displaystyle \lim_{s\to t}\dfrac{f(s)-f(t)}{g(s)-g(t)}, && t\notin D_{g}\cup C_{g},\\ 
             \vspace{3mm}
             \displaystyle\lim_{s\to t^{+}}\dfrac{f(s)-f(t)}{g(s)-g(t)}, && t\in D_{g}, \\
             \displaystyle\lim_{s\to b_n^+}\dfrac{f(s)-f(t)}{g(s)-g(t)}, && t\in (a_n,b_n)\subset C_{g},
         \end{array}
         \right.
     \end{equation}
     where $a_n$, $b_n$ are defined as in \eqref{cg}, provided that the corresponding limit exists. In such case, we say that $f$ is $g$-differentiable at $t$.
\end{definition}
Proposition \ref{tffemato} below will make it clear that Definition \ref{defderevery} is appropriate, but we need to recall some definitions first. Another natural notion in Stieltjes calculus is that of continuity with respect to the derivator $g$, see \cite{FrLo}. We say that a function $f:[a,b]\subset\mathbb{R}\longrightarrow\mathbb{R}$ is $g$-continuous at a point $t\in [a,b]$ if for every $\varepsilon>0$ there exists $\delta>0$ such that $|f(t)-f(s)|<\varepsilon$ whenever $s\in[a,b]$ and $|g(t)-g(s)|<\delta$.
    \noindent
    We say that $f$ is $g$-continuous on $[a,b]$, and denote it by $f\in\mathcal{C}_g\left([a,b]\right)$, if it is $g$-continuous at every $t\in[a,b]$. Surprisingly enough, $g$-continuous functions in compact intervals need not be bounded, see Example 3.3 in \cite{FrLo}, so we define
    $$\mathcal{BC}_g([a,b])=\{\mathcal{C}_{g}\left([a,b]\right) : f \ \text{is bounded}\},$$
    which is a Banach space (see \cite{FrLo}) with the sup norm
$$\|f\|_0=\sup_{t \in [a,b]}|f(t)|, \quad f \in \mathcal{BC}_g([a,b]).$$
Moreover, $g$-continuous functions in $[a,b]$ are regulated so, in particular, they are Borel measurable. Hence, $\mathcal{BC}_g([a,b]) \subset L^1_g([a,b])$. On the other hand, it is clear from Definition \ref{defacg} and Proposition \ref{bv} that $\mathcal{AC}_g([a,b]) \subset \mathcal{BC}_g([a,b])$.

It is important to note that we now have an expression for $f_{g}'$ at every point of $[a,b]$. Hence, we can take the $g$-derivatives of higher-order and consider the spaces of functions $k$-times continuously $g$-differentiable. Here we follow \cite{FeMaTo} again to define
$$
\mathcal{C}^{k}_g\left([a,b]\right)=\left\{f\in\mathcal{C}^{k-1}_g\left([a,b]\right)\ : \ (f_{g}^{(k-1)})'_{g}\in \mathcal{C}^{k-1}_g\left([a,b]\right) \right\},
$$
where $\mathcal{C}^{0}_g\left([a,b]\right)=\mathcal{C}_g\left([a,b]\right)$, $f_{g}^{(0)}=f$ and $f_{g}^{(k)}=(f_{g}^{(k-1)})'_{g}$.
For each $k\in\mathbb{N}$, we also define the set
$$\mathcal{BC}^{k}_{g}\left([a,b]\right)=\{f \in \mathcal{C}^{k}_{g}\left([a,b]\right) : f^{(n)}_g \in \mathcal{BC}_g([a,b]), \, n=0,1,2,\dots,k  \},$$
which is a Banach space (see \cite{FeMaTo}, Theorem 3.15) with the norm
$$
    \|f\|_k=\sum_{i=0}^k\|f^{(i)}_g\|_0, \quad f\in\mathcal{BC}^k_{g}\left([a,b]\right).
$$

\noindent
The space of $g$-absolutely continuous functions $\mathcal{AC}_{g}\left([a,b]\right)$ is a subset of $\mathcal{BC}_{g}\left([a,b]\right)$, see \cite[Proposition 5.5]{FrLo}. We also have $\mathcal{BC}^1_{g}\left([a,b]\right) \subset \mathcal{AC}_{g}\left([a,b]\right)$ along with some more information gathered in the following statement; see \cite[Lemma 3.13 and Lemma 3.14]{FeMaTo}.

\begin{proposition}
\label{tffemato}
Let $g:\mathbb{R}\longrightarrow\mathbb{R}$ be a derivator such that $a\notin N_{g}^{-}$ and $b\notin C_{g}\cup N_{g}^{+}\cup D_{g}$.

The embedding $\mathcal{BC}_g^1([a,b]) \hookrightarrow \mathcal{AC}_{g}([a,b])$ is continuous. Furthermore, for every $f \in \mathcal{BC}_g^1([a,b])$ we have
$$f(t)=f(a)+\int_{[a,t)}f'_g(s) \, dg(s) \quad \mbox{for all $t \in [a,b]$.}$$
Converserly, if $h \in \mathcal{BC}_g([a,b])$, then the function
$$H(t)=\int_{[a,t)}{h(s) \, dg(s)}, \quad t \in [a,b],$$
satisfies $H'_g(t)=h(t)$ for every $t \in [a,b]$ and, therefore, $H \in \mathcal{BC}^1_g([a,b])$.
\end{proposition}

Lastly, we include two results from \cite{FrLo} and \cite{frto} on the compactness of these function spaces, which we will need in the following sections.

\begin{theorem}\label{relcomp3}
    Let $\mathcal{S}\subset\mathcal{BC}_{g}\left([a,b]\right)$. If 
    \begin{enumerate}
       \item $S(t)$ is bounded for every $t\in[a,b]$ and 
       \item $S$ is uniformly $g$-equicontinuous, i.e., for every $\varepsilon>0$ and every $t\in[a,b]$ there exists $\delta>0$ such that 
        $$
        |f(t)-f(s)|<\varepsilon, \quad\text{for every}\,\,t,s\in[a,b]\,\,\text{such that}\,\,|g(t)-g(s)|<\delta,\,\,\text{and for every}\,\,f\in\mathcal{S},
        $$
        then $S$ is relatively compact.
    \end{enumerate}
    
\end{theorem}

\begin{theorem}\label{relcomp2}
Let $\mathcal{S}\subset\mathcal{BC}_{g}\left([a,b]\right)$. Then, $\mathcal{S}$ is a relatively compact subset of $(\mathcal{BC}_{g}\left([a,b]\right), \| \cdot \|_0)$ if it satisfies the following conditions:
    \begin{enumerate}
        \item $\mathcal{S}$ is $g$-equicontinuous, i.e., for every $\varepsilon>0$ and every $t\in[a,b]$ there exists $\delta>0$ such that 
        $$
        |f(t)-f(s)|<\varepsilon, \quad\text{for every}\,\,t,s\in[a,b]\,\,\text{such that}\,\,|g(t)-g(s)|<\delta,\,\,\text{and for every}\,\,f\in\mathcal{S}.
        $$
        \item The set $\{f(a):f\in\mathcal{S}\}$ is bounded.
        \item $\mathcal{S}$ has uniform right-hand side limits in $D_{g}\cap[a,b)$, i.e., for every $t_0\in D_{g}\cap[a,b)$ and every $\varepsilon>0$ there exists $\delta>0$ such that, for every $f\in\mathcal{S}$, we have
        $$
        |f(t)-f(t_0^+)|<\varepsilon, \quad \text{whenever}\,\, t\in(t_0,t_0+\delta).
        $$
        \item $\mathcal{S}$ has uniformly bounded jump discontinuities in in $D_{g}\cap[a,b)$, i.e., for every $t\in D_{g}\cap[a,b)$ there exists $\gamma_t>0$ such that, for every $f\in\mathcal{S}$, we have
        $$
        |f(t)-f(t^+)|<\gamma_t.
        $$
    \end{enumerate}
\end{theorem}
\begin{theorem}\label{relcomp}
    Let $\mathcal{S}\subset\mathcal{AC}_{g}\left([a,b]\right)$ be such that $\{f(a):f\in\mathcal{S}\}$ is bounded. Assume there exists a function $h\in L^{1}_{g}([a,b])$ such that
    $$
    |f'_{g}(t)|\le h(t), \quad \text{for $g$-almost all} \ t\in[a,b), \ \text{and for all}\ f\in\mathcal{S}.
    $$
    Then, $\mathcal{S}$ is a relatively compact subset of $(\mathcal{BC}_{g}\left([a,b]\right), \| \cdot \|_0)$.
\end{theorem}

Let us close this section by recalling that, for any derivator $g:[a,b]\to\mathbb{R}$, the Lebesgue--Stieltjes integration theory guarantees that 
\begin{equation}
\label{mvthm}
\int_{[a,b)}{f(s) \, dg(s)}=f(\sigma)(g(b)-g(a))
\end{equation}
for some $\sigma \in [a,b)$ provided that $f$ is uniformly continuous in $[a,b)$. This result is false if we only assume that $f$ is $g$--continuous, as the following simple example shows.

\begin{example}
Let $g(x)=x$ for $x \le 0$ and $g(x)=1+x$ for $x>0$. Clearly, $g$ is nondecreasing and left--continuous everywhere. We want to integrate with respect to $g$ the function $f(x)=1$ for $x \le 0$ and $f(x)=-1$ for $x>0$, which is $g$--continuous everywhere. We have
$$\int_{[-1,1)}{f(s) dg(s)}=\int_{-1}^0{ds}+\int_{\{0\}}{f(s)dg(s)}-\int_{0}^1{ds}=f(0)(g(0^+)-g(0))=1.$$
However, we have (\ref{mvthm}) for no $\sigma \in [-1,1]$, because
$$f(\sigma)(g(1)-g(-1))=3f(\sigma) \neq 1.$$
\end{example}
In this paper, we need the following mean value integral inequality.

\begin{proposition}
\label{mvii}
Let $g: \mathbb R \longrightarrow \mathbb R$ be a derivator and let $a,b \in \mathbb R$, $a<b$.

If $f \in L^1_g([a,b))$, then the sets
$$A_-=\left\{ \sigma \in [a,b) \, : \, \int_{[a,b)}{f(s) \, dg(s)} \le f(\sigma)(g(b)-g(a)) \right\}$$
and
$$A_+=\left\{ \sigma \in [a,b) \, : \, \int_{[a,b)}{f(s) \, dg(s)} \ge f(\sigma)(g(b)-g(a)) \right\}$$
have positive $g$--measure (in particular, they are not empty).
\end{proposition}

\noindent
{\bf Proof.} Let us prove that $\mu_g(A_-)>0$. The proof that $\mu_g(A_+)>0$ is analogous and we omit it. Firstly, observe that the result is trivial if $g(b)-g(a)=0$, so we assume that $g(b)-g(a)>0$. Secondly, the set $A_-$ is $g$--measurable because $f$ is a $g$--measurable function. Reasoning by contradiction, we assume that $\mu_g(A_-)=0$, i.e., 
$$\kappa=\dfrac{\int_{[a,b)}{f(s) \, dg(s)}}{g(b)-g(a)} > f(\sigma) \quad \mbox{for $g$--almost all $\sigma \in [a,b)$.}$$
By the strict monotonicity of the Lebesgue--Stieltjes integral on positive measure sets, we then have
$$\int_{[a,b)}{f(s) \, dg(s)} =\int_{[a,b)}{\kappa \, dg(\sigma)}>\int_{[a,b)}{f(\sigma) \, dg(\sigma)},$$
a contradiction.
\qed

\subsection{Strong and weak solutions}
The presence of two derivators in (\ref{eqg}), introduces some subtleties when it comes to defining the notion of a solution in the usual Carath\'eodory sense. 

Following \cite{FeMaPa}, where the case $g_1=g_2=Id$ is considered, we start with the analogous Carath\'eodory definition. Here and henceforth, we assume that $a\notin N_{g_1}^{-}$ and $b\notin C_{g_1}\cup N_{g_1}^{+}\cup D_{g_1}$, so that we can define derivatives with respect to $g_1$ at every point of $[a,b]$ by means of Definition \ref{defderevery}.

\begin{definition}\label{def:solg}
   A strong Carath\'eodory solution of the equation \eqref{eqg} is a function $x\in\mathcal{BC}^1_{g_1}(I)$ such that $(A\circ x)\cdot(\Phi\circ x'_{g_1})\in\mathcal{AC}_{g_2}(I)$ and
   $$
   ((A\circ x) \cdot(\Phi \circ x'_{g_1}))'_{g_2}(t) \in G(t,x(t),x'_{g_1}(t)) \,\ \text{for}\,\, g_2\text{-almost all}\,\, t \in I.
   $$
\end{definition}

Without additional assumptions, this definition quickly proves inadequate. In general, there is no hope of having strong solutions even if $G$ is constant and singlevalued.

\begin{example}
\label{exnosol}
Define $g_1(t)=t$ for $t \le 1/3$, $g_1(t)=1/3$ for $t \in [1/3,2/3]$, and $g_1(t)=t-1/3$ for $t>2/3$; also define $g_2(t)=t$ for all $t \in I$. The Stieltjes differential equation
\begin{equation}
\label{examplenosol}
\left(  x'_{g_1} \right)_{g_2}'(t)=1 \, \, \mbox{for $g_2$--a.a. $t \in I=[0,1]$,}
\end{equation}
has no strong solution. Indeed, assume in contrast that $x(t)$ is a solution in the sense of Definition \ref{def:solg}. Since $g_2$ is the identity, the corresponding Lebesgue--Stieltjes measure is the usual Lebesgue measure and then $x'_{g_1}$ is an absolutely continuous function on $I$ such that
$$\left(  x'_{g_1} \right)'(t)=1 \, \, \mbox{for a.a. $t \in I=[0,1]$.}$$
Therefore, 
\begin{equation}
\label{examplenosol2}
x'_{g_1}(t)=x'_{g_1}(0)+t \, \, \mbox{for all $t \in [0,1]$,}
\end{equation}
which does not belong to ${\cal BC}_{g_1}(I)$ because it is not constant in $[1/3,2/3]$.
\end{example}

Next, we introduce a more general definition, which also reduces to the definition given in \cite{FeMaPa} when $g_1=g_2=Id$.

\begin{definition}\label{def:solgw}
   A weak Carath\'edory solution of equation \eqref{eqg} is a function $x\in\mathcal{AC}_{g_1}(I)$ for which there exists $\xi  \in\mathcal{AC}_{g_2}(I)$ such that
   \begin{enumerate}
       \item[(1)] $\xi(t)=A(x(t))\Phi(x'_{g_1}(t))$ for $g_1$-almost all $t \in [a,b)$; and
       \item[(2)] $\xi'_{g_2}(t) \in G(t,x(t),\Phi^{-1}\left(\xi(t)/A(x(t)) \right)) \,\ \text{for}\,\, g_2\text{-almost all}\,\, t \in I.$
   \end{enumerate}
\end{definition}

Example (\ref{examplenosol}) has infinitely many weak solutions. Indeed, we know from (\ref{examplenosol2}) that any possible solution must satisfy
$$x_{g_1}'(t)=x'(t)=t+\lambda \quad \mbox{for a.a. $t \in (0,1/3)\cup (2/3,1)$, for some $\lambda \in \mathbb R$.}$$
Now it suffices to integrate with respect to $g_1$ to obtain the family of weak solutions given by
$$x(t)=\left\{
\begin{array}{cl}
t^2/2+\lambda \, t+\mu, & t \in [0,1/3],\\
1/18+\lambda/3+\mu, & t \in (1/3,2/3), \\
t^2/2+\lambda \, t+\mu -\lambda/3-1/6, & t \in [2/3,1], 
\end{array}
\right.$$
where $\lambda$ and $\mu$ are arbitrary real numbers.

\bigbreak

So far, we have seen that strong solutions can be rare while weak solutions can be too many. At this point, we decide to make an assumption which turns weak solutions into strong ones, namely,
${g_2}_{|I} \in {\cal C}_{g_1}(I)$.
\begin{proposition}
\label{prowest}
Assume that conditions (i), (ii), and (iii) are met. Let $x$ be a weak solution of (\ref{eqg}) and let $\xi$ be as in Definition \ref{def:solgw}.

If ${g_2}_{|I} \in {\cal C}_{g_1}(I)$ and there exists $d \in (0,c)$ such that
\begin{equation}
\label{acotacion}
\left| \dfrac{\xi(t)}{A \circ x(t)}\right| \le d \quad \mbox{for all $t \in I,$}
\end{equation}
then $x$ is a strong solution of (\ref{eqg}).
\end{proposition}

\noindent
{\bf Proof.} By definition, $\xi \in {\cal AC}_{g_2}(I) \subset {\cal BC}_{g_2}(I)$. The condition ${g_2}_{|I} \in {\cal C}_{g_1}(I)$ implies that $\xi \in {\cal BC}_{g_1}(I)$. We also have $A \circ x \in {\cal BC}_{g_1}(I)$. Moreover, there exists $m>0$ such that $A \circ x \ge m$ in $I$ and, as a result, $\xi / (A \circ x) \in {\cal BC}_{g_1}(I)$. The function $\Phi^{-1}$ is uniformly continuous in $[-d,d]$, hence, thanks to (\ref{acotacion}), we deduce that 
$$\Phi^{-1} \circ \dfrac{\xi} {A \circ x} \in  {\cal BC}_{g_1}(I).$$
Since $x \in {\cal AC}_{g_1}(I)$, for each $t \in [0,1]$ we have
$$x(t)=x(0)+\int_{[0,t)}{x'_{g_1}(s) \, dg_1(s)}=
x(0)+\int_{[0,t)}{\Phi^{-1}\left(\dfrac{\xi(s)}{A(x(s))}\right) \, dg_1(s)},$$
and then Proposition \ref{tffemato} ensures that
$$x'_{g_1}=\Phi^{-1} \circ \dfrac{\xi} {A \circ x} \quad 
\mbox{everywhere in $I$}.$$
Hence, $x \in {\cal BC}^1_{g_1}(I)$, $(A\circ x)\cdot(\Phi\circ x'_{g_1})=\xi \in {\cal AC}_{g_2}(I)$ and
$$((A\circ x)\cdot(\Phi\circ x'_{g_1}))'_{g_2}(t) \in G(t,x(t),x_{g_1}'(t)) \quad \mbox{for $g_2$--a.a. $t \in [a,b]$.}$$ \qed
Interestingly, condition (\ref{acotacion}) alone implies that weak solutions are strong in the case $g_1=g_2=Id$ considered in \cite{FeMaPa}, and this fact can be used to get a simpler proof of one of the main results in that reference, as we will do in our next section.

\section{An auxiliary problem}\label{sec3}
Following \cite{FeMaPa}, we consider the auxiliary Dirichlet problem
\begin{equation}\label{eqg2}
   ({\cal A}_{x} \cdot (\Phi \circ x'_{g_1}))'_{g_2}(t)\in F_{x}(t) \quad \mbox{for $g_2$-almost all $t \in I=[a,b]$,} \quad x(a)=v_a, \, \, x(b)=v_b,
\end{equation}
where
\begin{enumerate}
\item[$(i)$] $g_1, \, g_2:\mathbb R \longrightarrow \mathbb R$ are nondecreasing and left-continuous everywhere. Moreover, $a\notin N_{g_1}^{-}$ and $b\notin C_{g_1}\cup N_{g_1}^{+}\cup D_{g_1}$.
    \item[$(ii)$] $\Phi:\mathbb{R}\rightarrow(-c,c)$, $c>0$, is an increasing homeomorphism;
    \item[$(iii)$] ${\cal A}:\mathcal{AC}_{g_1}(I)\rightarrow\mathcal{BC}_{g_1}(I)$ is a uniformly continuous operator which maps bounded sets into relatively compact sets and there exist $m,M>0$ such that
    \begin{equation}\label{cotA}
    0<m\leq A_{x}(t)\leq M \quad \text{for every}\, \, x \in \mathcal{AC}_{g_1}(I)\,\, \text{and}\,\, t\in I;
    \end{equation}
    \item[$(iv)$] $F:\mathcal{AC}_{g_1}(I)\rightarrow\mathcal{P}(L^{1}_{g_2}(I))$ is a w-u.s.c. multivalued mapping with non-empty, closed and convex values, and there exists a function $\psi \in L^{1}_{g_2}(I)$ such that for every $x\in\mathcal{AC}_{g_1}(I)$ we have
    \begin{equation}\label{cotaF}
        \| F_{x}(t)\|=\sup\{|f(t)| : f \in F_x\}\leq\psi(t) \quad\mbox{for $g_2$-almost all $t\in I$.}
    \end{equation}
\end{enumerate} 
\bigbreak
Consistently with Definition \ref{def:solg}, we define a strong solution of problem \eqref{eqg2} as a function $x\in\mathcal{BC}^1_{g_1}(I)$ such that $x(a)=v_a$, $x(b)= v_b$, ${\cal A}_{x} \cdot (\Phi \circ x'_{g_1})\in\mathcal{AC}_{g_2}(I)$ and there exists $f \in F_x$ such that $$({\cal A}_{x} \cdot (\Phi \circ x'_{g_1}))'_{g_2}(t)=f(t) \quad \mbox{for $g_2$-a.a. $t\in I$.}$$
Similarly, a weak solution of problem \eqref{eqg2} is a function $x\in\mathcal{AC}_{g_1}(I)$ such that $x(a)=v_a$, $x(b)= v_b$, and there exists $\xi \in \mathcal{AC}_{g_2}(I)$ such that $\xi={\cal A}_{x} \cdot (\Phi \circ x'_{g_1})$ for $g_1$-almost all $t \in I$ and for some $f \in F_x$ we have $$\xi'_{g_2}(t)=f(t) \quad \mbox{for $g_2$-a.a. $t\in I$.}$$
\bigbreak

In these conditions, we can integrate with respect to $g_2$ to obtain
\begin{equation}
    \label{xi}
\xi(t)=\xi(a)+\int_{[a,t)}{f(s) \, dg_2(s)} \quad \mbox{for all $t \in I$},
\end{equation}
and, therefore, for $g_1$-almost all $t \in I$ we have
$$x'_{g_1}(t)=\Phi^{-1} \left( \dfrac{\xi(t)}{{\cal A}_x(t)}\right)= \Phi^{-1} \left( \dfrac{\xi(a)+\int_{[a,t)}{f(s) \, dg_2(s)}}{{\cal A}_{x}(t)} \right).$$
Integrating with respect to $g_1$, we get
\begin{equation}
\label{pf}
 x(t)=v_a+\int_{[a,t)}{\Phi^{-1} \left( \dfrac{\xi(a)+\int_{[a,r)}{f(s) \, dg_2(s)}}{{\cal A}_{x}(r)} \right) \, dg_1(r)}  \quad \mbox{for all $t \in I$,}
\end{equation}
and the boundary condition $x(b)=v_b$ implies that $\tau=\xi(a)$ is implicitly determined by the algebraic equation
\begin{equation}
\label{impint}
\int_{[a,b)}{\Phi^{-1} \left( \dfrac{\tau+\int_{[a,r)}{f(s) \, dg_2(s)}}{{\cal A}_{x}(r)} \right) \, dg_1(r)}=v_b-v_a.
\end{equation}
Conversely, if for some $x \in {\cal AC}_{g_1}(I)$ there exists $f \in F_x$ such that $x$ satisfies (\ref{pf}) for $\xi(a)=\tau$, and $\tau$ given by (\ref{impint}), then $x$ is a weak solution of (\ref{eqg2}), with the corresponding function $\xi(t)$ defined by means of (\ref{xi}).

Summing up, weak solutions of (\ref{eqg2}) can be seen as fixed points of a multivalued integral operator defined by means of (\ref{pf}) and (\ref{impint}). However, such a multivalued operator does not assume convex values and that is why we need Theorem \ref{teo:fp1} to study it. 
\bigbreak
In the following theorem we prove the existence of at least one solution for problem \eqref{eqg2} and the $g_1$-slope
$$
s^*=\dfrac{v_b-v_a}{g_1(b)-g_1(a)}
$$
plays a fundamental role. Note that the second part of condition $(i)$ implies that $g_1(b)-g_1(a) \neq 0$.

\begin{theorem}\label{teo:fun}
    Let $g_1$, $g_2$, $\Phi$, ${\cal A}$ and $F$ satisfy conditions $(i)$-$(iv)$ at the beginning of Section 3. 
    
    If
    \begin{equation}\label{ms}
        M|\Phi(s^*)| <cm-2\|\psi\|_{L^{1}_{g_2}},
    \end{equation}
    then problem \eqref{eqg2} has at least one weak solution, which is a strong solution provided that ${g_2}_{|I} \in {\cal C}_{g_1}(I)$.
    \begin{proof}
        Let $\widehat{F}:\mathcal{AC}_{g_1}(I)\rightarrow\mathcal{P}(\mathcal{BC}_{g_2}(I))$ be the Carathéodory operator corresponding to $F$, that is,
        $$
        \widehat{F}_{x}=\left\lbrace y \in \mathcal{BC}_{g_2}(I)\,\,  : \,\,  y(t)=\int_{[a,t)}f(s)\,dg_2(s),\, t\in I, \,\, \text{for some}\,\, f\in F_{x}\right\rbrace,\quad x \in \mathcal{AC}_{g_1}(I).
        $$
        Since $F$ assumes nonempty convex values, so does $\widehat{F}$.
        
        \noindent
        {\it Claim 1: $\widehat{F}$ is u.s.c.} To prove it, we use Proposition \ref{prohupa}. First, we prove that $\widehat{F}$ assumes closed values. We fix $x\in\mathcal{AC}_{g_1}(I)$, take a sequence $\{y_n\}_{n\in\mathbb{N}}$ in $\widehat{F}_x$ such that $y_n\to y$ uniformly in $I$, and we have to prove that $y\in\widehat{F}_x$. Note that for each $n\in\mathbb{N}$, there exists $f_n\in F_x$ such that $y_n(t)=\int_{[a,t)}f_n(s)\,dg_2(s)$, $t\in I$. Since $\{f_n\}_{n\in\mathbb{N}}\subset F_x$, condition \eqref{cotaF} implies that it is a uniformly integrable sequence, so by the Dunford-Pettis Theorem it has a subsequence $\{f_{n_k}\}_{k\in\mathbb{N}}$ which is weakly convergent to a function $f\in L_{g_2}^{1}(I)$. Therefore, for each fixed $t \in I$, the mapping $\chi_{[a,t)}\in L^\infty_{g_2}(I)$, and we have
        $$
        y(t)=\lim_{k \to \infty}y_{n_k}(t)=\lim_{k\to\infty}\int_{[a,b)}f_{n_k}(s)\chi_{[a,t)}(s)\,dg_2(s)=\int_{[a,t)}f(s)\,dg_2(s).
        $$
        By Mazur's Lemma, there is a sequence $\left\{g_m\right\}_{m \in \mathbb N}$ made up of convex linear combinations of the $f_{n_k}$'s that converges strongly to $f$ in $L^1_{g_2}(I)$, and this implies that $f \in F_x$ because $F_x$ is convex and closed in $L^1_{g_2}(I)$. Hence, $y \in \hat{F}_x$.

        \noindent
        Second, we need to prove that $\widehat{F}$ is locally compact, but much more than that is true. Indeed, we shall prove that $\widehat{F}(\mathcal{AC}_{g_1}(I))$ is relatively compact in $\mathcal{BC}_{g_2}(I)$. Clearly, $\widehat{F}(\mathcal{AC}_{g_1}(I)) \subset \mathcal{AC}_{g_2}(I)$, and for each $y \in \widehat{F}(\mathcal{AC}_{g_1}(I))$, we have $y(a)=0$ and $|y'_{g_2}(t)| \le \psi(t)$ for $g_2$-a.a. $t \in [a,b)$. Hence, Theorem \ref{relcomp} guaranties that $\widehat{F}(\mathcal{AC}_{g_1}(I))$ is relatively compact in $\mathcal{BC}_{g_2}(I)$.

        \noindent
        Third and last, we show that $\widehat{F}$ has a closed graph. Consider the graph
        $$
        \mathrm{Gr}(\widehat{F})=\left\lbrace(x,y)\in \mathcal{AC}_{g_1}(I) \times \mathcal{BC}_{g_2}(I)\, : \,y\in \widehat{F}_x\right\rbrace,
        $$
        and let $\{(x_n,y_n)\}_{n\in\mathbb{N}}$ be a sequence in $\mathrm{Gr}(\widehat{F})$ such that $(x_n,y_n)\to(x,y)$. We need to prove that $y\in \widehat{F}_x$. By the definition of $\widehat{F}$, for each $n\in\mathbb{N}$ there exists $f_n\in F_{x_n}$ such that
        $y_n(t)=\int_{[a,t)}f_{n}(s)\,dg_2(s)$. Again using condition \eqref{cotaF} and the Dunford-Pettis Theorem, we deduce that the sequence $\{f_n\}_n$ has a subsequence $\{f_{n_k}\}_k$ that converges weakly to a function $f\in L_{g_2}^{1}(I)$ which belongs to $F_x$ because $F$ is w-u.s.c. Hence, for each fixed $t \in I$, we have
        $$
        y(t)=\lim_{k\to\infty}y_{n_k}(t)=\lim_{k\to\infty}\int_{[a,b)}f_{n_k}(s)\chi_{[a,t)}(s)\,dg_2(s)=\int_{[a,t)}f(s)\,dg_2(s),
        $$
    thus proving that $y\in{\widehat{F}_x}$ as desired. 
   
\bigbreak
        Denote now
        $$
        \Delta=\left\lbrace y \in \mathcal{BC}_{g_2}(I) \,\, :\,\, \|y\|_{0}\leq\|\psi\|_{L^1_{g_2}}\right\rbrace.
        $$
        \textit{Claim 2:} For each $x\in\mathcal{AC}_{g_1}(I)$ and each $y\in \Delta$ there exists a unique $\beta_{xy}\in\mathbb{R}$ such that
        \begin{equation}\label{betaxy}
            \int_{[a,b)}\Phi^{-1}\left(\dfrac{\beta_{xy}+y(t)}{{\cal A}_{x}(t)}\right)\,dg_1(t)=v_b-v_a.
        \end{equation}
        \bigbreak
        To prove this, fix $x\in\mathcal{AC}_{g_1}(I)$ and $y\in \Delta$, and consider the function
        \begin{equation*}
        \varphi(\tau)=\int_{[a,b)}\Phi^{-1}\left(\dfrac{\tau+y(t)}{{\cal A}_{x}(t)}\right)\,dg_1(t),\quad \tau\in J=\left(-cm+\|\psi\|_{L^1_{g_2}},cm-\|\psi\|_{L^1_{g_2}}\right).
        \end{equation*}

        \noindent
        Note that for each fixed $\tau_0 \in J$ we can specify a neighborhood $V_0$ and a constant $\rho \in (0,cm-\|\psi\|_{L^1_{g_2}})$ such that $\tau \in V_0$ implies $|\tau| \le \rho$, and therefore
        \begin{equation*}
        \left|\dfrac{\tau+y(t)}{{\cal A}_{x}(t)}\right|\le \dfrac{\rho+\|\psi\|_{L^1_{g_2}}}{m}<
        \dfrac{cm-\|\psi\|_{L^1_{g_2}}+\|\psi\|_{L^1_{g_2}}}{m}= c, \, \, t \in I.
        \end{equation*}
        Hence, $\varphi(\tau)$ is well defined for all $\tau \in V_0$ because the function under the integral sign is Borel measurable (because it is regulated) and bounded, hence Lebesgue--Stieltjes integrable with respect to $g_1$. Moreover, the Dominated Convergence Theorem can be used to prove that $\varphi$ is continuous at $\tau_0$. Obviously, $\varphi$ is strictly increasing.

        \noindent
        Now, by condition \eqref{ms}, for every $t\in I$ we have 
        \begin{equation*}
            \dfrac{-cm+\|\psi\|_{L^1_{g_2}}+y(t)}{A_{x}(t)}\leq \dfrac{-cm+2\|\psi\|_{L^1_{g_2}}}{M}<\Phi(s^*)<\dfrac{cm-2\|\psi\|_{L^1_{g_2}}}{M}\leq\dfrac{cm-\|\psi\|_{L^1_{g_2}}+y(t)}{A_{x}(t)}.
        \end{equation*}
        For a sufficiently small $r>0$, we still have
 \begin{equation*}
             \dfrac{r-cm+2\|\psi\|_{L^1_{g_2}}}{m}<\Phi(s^*)<\dfrac{cm-r-2\|\psi\|_{L^1_{g_2}}}{M},
        \end{equation*}
which implies that for any $t$ we have
         \begin{equation*}
            -c<\dfrac{r-cm+\|\psi\|_{L^1_{g_2}}+y(t)}{{\cal A}_{x}(t)} <\Phi(s^*)< \dfrac{cm-r-\|\psi\|_{L^1_{g_2}}+y(t)}{{\cal A}_{x}(t)}<c,
        \end{equation*}
        so, by the monotonicity of $\Phi^{-1}$,
        \begin{equation*}
            \Phi^{-1}\left(\dfrac{r-cm+\|\psi\|_{L^1_{g_2}}+y(t)}{{\cal A}_{x}(t)}\right)<s^*<\Phi^{-1}\left(\dfrac{cm-r-\|\psi\|_{L^1_{g_2}}+y(t)}{{\cal A}_{x}(t)}\right),
        \end{equation*}
        and, subsequently,
        \begin{equation*}
           \varphi\left(r-cm+\|\psi\|_{L^1_{g_2}}\right) <v_b-v_a<\varphi\left(cm-r-\|\psi\|_{L^1_{g_2}}\right).
        \end{equation*}
        Since $\varphi$ is strictly increasing and continuous, we deduce the existence of a unique $$\beta_{xy}\in [r-cm+\|\psi\|_{L^1_{g_2}},cm-r-\|\psi\|_{L^1_{g_2}}] \subset J$$ such that $\varphi(\beta_{xy})=v_b-v_a$, so Claim $1$ is proved.
\bigbreak

        Next, we define the singlevalued operator $$H:(x,y) \in \mathcal{AC}_{g_1}(I)\times\Delta\longmapsto H_{xy} \in \mathcal{AC}_{g_1}(I),$$ where
        \begin{equation}\label{hxy}
        H_{xy}(t)=v_a+
            \int_{[a,t)}\Phi^{-1}\left(\dfrac{\beta_{xy}+y(s)}{{\cal A}_{x}(s)}\right)\,dg_1(s), \, \, t \in I.
        \end{equation}
         Let us prove that $H$ is well defined. Proposition \ref{mvii} applied to (\ref{betaxy}) implies the existence of points $\sigma_1, \sigma_2 \in I$ such that
        $$
        \Phi^{-1}\left(\dfrac{\beta_{xy}+y(\sigma_1)}{{\cal A}_{x}(\sigma_1)}\right) \le \dfrac{v_b-v_a}{g_1(b)-g_1(a)}=s^* \le \Phi^{-1}\left(\dfrac{\beta_{xy}+y(\sigma_2)}{{\cal A}_{x}(\sigma_2)}\right).
        $$
      As a result, we get
        $$
        \beta_{xy}\leq {\cal A}_{x}(\sigma_1)\Phi(s^*) -  y(\sigma_1) \leq M|\Phi(s^*)|+\|\psi\|_{L_{g_2}^1},
        $$
        and
 $$
        -\beta_{xy}\leq -{\cal A}_{x}(\sigma_2)\Phi(s^*) -  y(\sigma_2) \leq M|\Phi(s^*)|+\|\psi\|_{L_{g_2}^1},
$$and therefore
\begin{equation}
\label{betab}
 |\beta_{xy}| \leq M|\Phi(s^*)|+\|\psi\|_{L_{g_2}^1}.
\end{equation}
        Hence, for each $x\in\mathcal{AC}_{g_1}(I)$ and each $y\in\Delta$, we have
        \begin{equation}
        \label{et7}
        |\beta_{xy}|+|y(t)| \le M |\Phi(s^*)|+2\|\psi\|_{L_{g_2}^1}\quad \mbox{for all $t \in I$,}
        \end{equation}
       \begin{equation}
       \label{etiq}
        \left|\dfrac{\beta_{xy}+y(t)}{{\cal A}_{x}(t)}\right|\leq \dfrac{M|\Phi(s^*)|+2\|\psi\|_{L_{g_2}^1}}{m}:=d \quad \mbox{for all $t \in I$,}
        \end{equation}
        and, by condition \eqref{ms}, $d<c$. Therefore, the function $\phi(t)=\Phi^{-1}(({\beta_{xy}+y(t)})/{{\cal A}_{x}(t)})$ is bounded and regulated, so it is $g_1$-integrable on $I$. Hence, $H_{xy}\in {\cal AC}_{g_1}(I)$. 
        
        Moreover, since $\Phi^{-1}$ is bounded on $[-d,d]$, we deduce the existence of a constant $k>0$ such that $\|H_{xy}\|_{{\cal AC}_{g_1}}\leq k$ for every $x\in\mathcal{AC}_{g_1}(I)$ and every $y\in\Delta$. In other words,
        $$
        H(\mathcal{AC}_{g_1}(I)\times\Delta)\subset \overline{B(0;k)}_{\mathcal{AC}_{g_1}}.
        $$
        \textit{Claim 3:} $H$ is uniformly continuous on $\mathcal{AC}_{g_1}(I)\times\Delta$.

        \noindent
        Let $x_1,x_2\in\mathcal{AC}_{g_1}(I)$, $y_1,y_2\in\Delta$, and denote $\beta_1=\beta_{x_1y_1}$, $\beta_2=\beta_{x_2y_2}$. By applying Proposition \ref{mvii} to the integral
        $$
        \int_{[a,b)}\left(\Phi^{-1}\left(\dfrac{\beta_{1}+y_1(t)}{{\cal A}_{x_1}(t)}\right)-\Phi^{-1}\left(\dfrac{\beta_{2}+y_2(t)}{{\cal A}_{x_2}(t)}\right)\right)\,dg_1(t)=0,
        $$
        we deduce the existence of a point $t^*\in I$ such that
        $$
        \Phi^{-1}\left(\dfrac{\beta_{1}+y_1(t^*)}{{\cal A}_{x_1}(t^*)}\right) \le \Phi^{-1}\left(\dfrac{\beta_{2}+y_2(t^*)}{{\cal A}_{x_2}(t^*)}\right).
        $$
        Hence, since $\Phi$ is increasing,
        $$
        {\cal A}_{x_2}(t^*)\beta_1+{\cal A}_{x_2}(t^*)y_1(t^*) \le {\cal A}_{x_1}(t^*)\beta_2+{\cal A}_{x_1}(t^*)y_2(t^*),
        $$
        so
        $$
        m(\beta_1-\beta_2) \le {\cal A}_{x_2}(t^*)(\beta_1-\beta_2) \le (\beta_2+y_1(t^*))\left({\cal A}_{x_1}(t^*)- {\cal A}_{x_2}(t^*)\right)+{\cal A}_{x_1}(t^*)(y_2(t^*)-y_1(t^*)).
        $$
        As a consequence, we deduce from (\ref{et7}) and (\ref{ms}) that
        \begin{equation}\label{beta12}
            m|\beta_1-\beta_2|\leq cm\|{\cal A}_{x_1}-{\cal A}_{x_2}\|_{0}+M\|y_2-y_1\|_{0}.
        \end{equation}
        Moreover, for every $t\in I$ we have that
        \begin{equation}\label{beta12_2}
        \left|\dfrac{\beta_{1}+y_1(t)}{{\cal A}_{x_1}(t)}-\dfrac{\beta_{2}+y_2(t)}{{\cal A}_{x_2}(t)}\right|\leq\dfrac{1}{m^2}\left\lbrace M\left(|\beta_1-\beta_2|+\|y_1-y_2\|_0\right)+cm\|{\cal A}_{x_1}-{\cal A}_{x_2}\|_{0}\right\rbrace.
        \end{equation}
        Now, since $\Phi^{-1}$ is uniformly continuous on $[-d,d]$, for any fixed $\varepsilon>0$ there exists a $\delta>0$ such that
        \begin{equation}\label{Phicont}
        |\Phi^{-1}(s_1)-\Phi^{-1}(s_2)|<\min\left\{\dfrac{\varepsilon}{2},\dfrac{\varepsilon}{2(g_1(b)-g_1(a))}\right\}, 
        \end{equation}
        for every $s_1,s_2\in [-d,d]$ such that $|s_1-s_2|<\delta$. Since $\mathcal{A}$ is uniformly continuous, we can find $\eta>0$ such that
        \begin{equation}\label{Acont}
        \|{\cal A}_{x_2}-{\cal A}_{x_1}\|_{0}<\dfrac{\delta m^2}{2c(m+M)},\quad \text{if}\ \ \|x_2-x_1\|_{\mathcal{AC}_{g_1}}<\eta.
        \end{equation}
        Now, take $\rho=\frac{\delta m^3}{2M(m+M)}$. Inequalities \eqref{beta12}, \eqref{beta12_2} and \eqref{Phicont} imply that
        $$
        \max_{t\in I}\left|\dfrac{\beta_{1}+y_1(t)}{{\cal A}_{x_1}(t)}-\dfrac{\beta_{2}+y_2(t)}{{\cal A}_{x_2}(t)}\right|\leq\dfrac{c(M+m)\|{\cal A}_{x_1}-{\cal A}_{x_2}\|_{0}+M(M+m)/m\|y_1-y_2\|_0}{m^2}<\delta,
        $$
        whenever $\|x_1-x_2\|_{\mathcal{AC}_{g_1}}<\eta $ and $\|y_1-y_2\|_0<\rho$. Hence
        $$
        \|H_{x_1y_1}-H_{x_2y_2}\|_{\mathcal{AC}_{g_1}}\leq\int_{[a,b)}\left|\Phi^{-1}\left(\dfrac{\beta_{1}+y_1(s)}{{\cal A}_{x_1}(s)}\right)-\Phi^{-1}\left(\dfrac{\beta_{2}+y_2(s)}{{\cal A}_{x_2}(s)}\right)\right|\, dg_1(s)<\varepsilon.
        $$
        so $H$ is uniformly continuous.

        Consider now the convex subset of $\mathcal{BC}_{g_2}(I)$ defined by
        $$
        D=\mathrm{co}\, \hat{F}(\mathcal{AC}_{g_1}(I)),
        $$
        which, by condition \eqref{cotaF}, is included in $\Delta$.

\bigbreak
        \noindent
        \textit{Claim 4:} $H\left(\overline{B(0;k)}_{\mathcal{AC}_{g_1}}\times D\right)$ is relatively compact in $\mathcal{AC}_{g_1}(I)$.\\
        Let $\{h_n\}_{n \in \mathbb N}$ be an arbitrary sequence in $ H\left(\overline{B(0;k)}_{\mathcal{AC}_{g_1}}\times D\right)$; we have to prove that some subsequence is convergent in $\mathcal{AC}_{g_1}(I)$. To prove it, note that for each $n \in \mathbb N$ there exist $x_n \in\overline{B(0;k)}_{\mathcal{AC}_{g_1}}$ and $y_n\in D$, such that
       \begin{equation}
           \label{h}
        h_n(t)=v_a+
            \int_{[a,t)}\Phi^{-1}\left(\frac{\beta_{n}+y_n(s)}{{\cal A}_{x_n}(s)}\right)\,dg_1(s), \, \, t \in I,
      \end{equation}
      where $\beta_n=\beta_{x_ny_n}$ is as in Claim 2. Obviously, $h_n(a)=v_a$ for all $n \in \mathbb N$. Moreover, 
      \begin{equation}
      \label{hprimcot}
        (h_n)'_{g_1}(t)=\Phi^{-1}\left(\frac{\beta_{n}+y_n(t)}{{\cal A}_{x_n}(t)}\right) \quad \mbox{for $g_1$-a.a. $t \in I$ and all $n \in \mathbb N$ }
        \end{equation}
        and, since $\Phi^{-1}$ is continuous on $[-d,d]$, with $d$ as in (\ref{etiq}), there exists $R>0$ such that $|(h_n)'_{g_1}(t)|\leq R$ for $g_1$-a.a. $t\in I$. Since this bound does not depend on the particular ${h_n}$ selected, and every constant function belongs to $L_{g_1}^{1}(I)$, we deduce from Theorem \ref{relcomp} that the sequence $\{h_n\}_{n \in \mathbb N}$ is relatively compact in $\mathcal{BC}_{g_1}(I)$. As a result, there exists $\gamma \in \mathcal{BC}_{g_1}(I)$ and a subsequence of $\{h_n\}_{n \in \mathbb N}$, which we denote again by $\{h_n\}_{n \in \mathbb N}$, such that $h_n \to \gamma$ uniformly in $I$. 
        
        Let us prove that, in fact, $\gamma \in \mathcal{AC}_{g_1}(I)$ and a subsequence of $\{h_n\}_{n \in \mathbb N}$ tends to $\gamma$ in the $\mathcal{AC}_{g_1}$-norm.

        For each $y_n\in D=\mbox{co} \widehat{F}(\mathcal{AC}_{g_1}(I))$ there exists $m(n) \in \mathbb N$, $\lambda_{n,i}\in[0,1]$, $u_{n,i}\in \widehat{F}(\mathcal{AC}^1_{g_1}(I))$ for $i\in\{1,\dots,m(n)\}$ such that
        \begin{equation}
            y_n=\sum_{i=1}^{m(n)}\lambda_{n,i}\,u_{n,i}.
        \end{equation}
       Hence, $y_n(a)=0$ for all $n \in \mathbb N$, and from the definition of $\widehat{F}$ and the condition \eqref{cotaF}, we get
        \begin{equation*}
            |(y_n)'_{g_2}(t)|\leq\sum_{i=1}^{m(n)}\lambda_{n,i}|(u_{n,i})'_{g_2}(t)|\leq \psi(t) \quad \mbox{for $g_2$-a.a. $t \in I$.}
        \end{equation*}
        Therefore, Theorem \ref{relcomp} guaranties that the sequence $\{y_n\}_{n \in \mathbb N}$
        is relatively compact in $\mathcal{BC}_{g_2}(I)$, so there is a subsequence $\{y_{n_k}\}_k$ and $\xi \in \mathcal{BC}_{g_2}(I)$, such that $\{y_{n_k}\}_k$ tends to $\xi$ uniformly in $I$.

        By condition $(iii)$, the sequence $\{A \circ x_{n_k}\}_k$ is relatively compact in $\mathcal{BC}_{g_1}(I)$, so it has a subsequence, which we denote again by $\{A \circ x_{n_k}\}_k$, which converges uniformly in $I$ to certain $\eta \in \mathcal{BC}_{g_1}(I).$ Moreover, $A(x_{n_k}(t))\ge m$ for all $t \in I$ and all $k \in \mathbb N$.

        Finally, we may assume w.l.o.g., that $\beta_{n_k} \to \beta_{\infty} \in \mathbb R$, and therefore, for $g_1$-a.a. $t \in I$, we have
        $$\lim_{k \to \infty} (h_{n_k})'_{g_1}(t)=\lim_{k\to \infty} \Phi^{-1} \left(\frac{\beta_{n_k}+y_{n_k}(t)}{A(x_{n_k}(t))}\right)=\Phi^{-1} \left(\frac{\beta_{\infty}+\xi(t)}{\eta(t)}\right)\equiv \varphi(t).$$
        Observe that $\varphi(t)$ is well--defined for every $t \in I$, and $\varphi$ is Borel-measurable because it is the pointwise limit of a sequence of regulated functions (hence Borel-measurable functions). Moreover, (\ref{hprimcot}) ensures that $\{(h_{n_k})'_{g_1}\}_k$ tends to $\varphi$ in $L^1_{g_1}(I)$, so for each $t \in I$ we have
        $$\gamma(t)=\lim_{k \to \infty}h_{n_k}(t)=\lim_{k \to \infty}\left( v_a+ \int_{[a,t)}\Phi^{-1}
        \left( \frac{\beta_{n_k}+y_{n_k}(t)}{A(x_{n_k}(t))}\right) dg_1(s)\right)=v_a +\int_{[a,t)}{\varphi(s) \, dg_1(s)}.$$
    As a result, $\gamma \in \mathcal{AC}_{g_1}(I)$ and
$h_{n_k} \to \gamma$ in $\mathcal{AC}_{g_1}(I)$. The proof of Claim 4 is complete.       
        
\bigbreak

        Finally, consider the convex and compact subset of $\mathcal{AC}_{g_1}(I)$ given by
        $$
        S=\overline{\mathrm{co}} \, \overline{H\left(\overline{B(0;k)}_{\mathcal{AC}_{g_1}}\times D\right)},
        $$
        and the operator $\Gamma:S\rightarrow\mathcal{P}(S)$, defined as
        $$
        \Gamma_{x}=H(x,\hat{F}_x)=\left\{H_{xy}\, : \,   y\in\hat{F}_x\right\}.
        $$
        We will check that $\Gamma$ has a closed graph. Let $\{(x_n,z_n)\}_n\subset\mathrm{Gr}(\Gamma)$ a sequence such that $||x_n-x||_{\mathcal{AC}_{g_1}}\to0$ and $||z_n-z||_{\mathcal{AC}_{g_1}}\to0$, and let us see that $(x,z)\in\mathrm{Gr}(\Gamma)$. By definition of $\Gamma$, for each $n\in\mathbb{N}$ there exists $y_n\in\hat{F}_{x_n}\subset\mathcal{AC}_{g_2}(I)$ such that $z_n=H(x_n,y_n)$. Since $\{y_n(a):n\in\mathbb{N}\}=\{0\}$ is bounded and $|(y_n)'_{g_2}(t)|\leq\psi(t)$ for $g_2$-almost all $t\in I$, $\{y_n\}_n$ is relatively compact in $\mathcal{BC}_{g_2}(I)$, so we deduce the existence of a function $y\in\mathcal{BC}_{g_2}(I)$ such that $\|y_n-y\|_{0}\to0$, passing to a subsequence if necessary. Since $\widehat{F}$ is u.s.c., $y\in\widehat{F}_x$. Moreover, since $H$ is uniformly continuous and $(x_n,y_n)\to (x,y)$, then $z_n=H(x_n,y_n)\to H(x,y)$, so $(x,z)\in \mathrm{Gr}(\Gamma)$. Hence, $\Gamma$ has a closed graph, so we can apply Proposition \ref{teo:fp1} to establish the existence of a fixed point for $\Gamma$, $x\in\Gamma_x$. Lastly, we deduce from equations \eqref{betaxy} and \eqref{hxy} that $x$ is a weak solution for problem \eqref{eqg2}, which turns into a strong solution if ${g_2}_{|I} \in {\cal C}_{g_1}(I)$ thanks to the inequality (\ref{etiq}) and the same arguments in the proof of Proposition \ref{prowest}.
     \end{proof}
\end{theorem}
\begin{example}
In order to give a proper motivation to the upcoming example, we should start by recalling the classical damped pendulum equation
\begin{equation}
\label{pend}
x''=k \sin (x)-\mu \, x',
\end{equation}
where $k$ and $\mu$ are positive. If the pendulum is attached to a rough wall, then it is reasonable to replace the usual linear damping term with Coulomb's dry friction and the differential equation with a differential inclusion, namely
\begin{equation}
\label{pend2}
x'' \in k \sin (x)-\mu \, \mbox{Sgn} \,(x'),
\end{equation}
where 
$$\mbox{Sgn} \,(z) = 
\left\{
\begin{array}{cc}
\{-1\},   & \mbox{if $z<0$,} \\
\left[-1,1\right],   & \mbox{if $z=0$,} \\
\{1\},   & \mbox{if $z>0$.} 
\end{array}
\right.
$$
Our theory applies even to more general versions of (\ref{pend2}) which can be useful from a physical viewpoint. Consider, for instance, the $\Phi$-laplacian version
\begin{equation}
\label{pend3}
(\Phi \circ x')' \in k \sin (x)-\mu \, \mbox{Sgn} \,(x'),
\end{equation}
where $\Phi$ is as in $(ii)$. Under suitable conditions, see \cite{capo}, we have
$$(\Phi \circ x')'=\Phi'(x')x'',$$
and we could rewrite the inclusion (\ref{pend3}) as
\begin{equation}
\label{pend4}x'' \in \dfrac{1}{\Phi'(x')} \left(k \sin (x)-\mu \, \mbox{Sgn} \,(x')\right),\end{equation}
where the term $1/\Phi'(x')$, which does not need to be continuous, can be seen as a speed dependent external control on the total forces. The main difficulty of (\ref{pend4}) with respect to (\ref{pend3}) is that the right-hand side of (\ref{pend4}) is not necessarily u.s.c. with respect to $x'$.

Last but not least, we introduce in (\ref{pend3}) a Stieltjes derivative
\begin{equation}
\label{pend5}
(\Phi \circ x')_g' \in k \sin (x)-\mu \, \mbox{Sgn} \,(x'), \quad t \in [0,1],
\end{equation}
where $g(t)=t+H(t-t_0)$, where $t_0 \in (0,1)$ is fixed and $H$ is the Heaviside function, namely $H(s)=0$ for $s \le 0$ and $H(s)=1$ for $s>0$. In this case, the $g$-derivative is the usual one for every $t \neq t_0$, but at $t=t_0$ the function $x'$ has a jump discontinuity. This is useful if we want to model the motion of the pendulum which is suddenly hit at time $t_0$, and thus changes its speed instantly and discontinuously.

Theorem \ref{teo:fun} can be applied to prove that the inclusion (\ref{pend5}) with Dirichlet boundary conditions has at least one weak solution. Specifically, if we consider \eqref{pend5} equipped with $x(0)=0$ and $x(1)=1$, then \eqref{ms} holds if $|\Phi(1)|<c-2(k+\mu)$, so the problem has at least one weak solution. However, there is not necessarily a strong solution, because $g$ is not continuous with respect to the identity function, the other derivator involved in the left-hand side of (\ref{pend5}).

\end{example}

\section{Existence of solution for the Dirichlet problem}\label{sec4}
In this section, we consider the Dirichlet problem
\begin{equation}\label{DP}
   ((A \circ x) \cdot (\Phi \circ x'_{g_1}))'_{g_2}(t)\in G(t,x(t),x'_{g_1}(t)), \quad t \in I=[a,b], \quad x(a)=v_a,\, \, x(b)=v_b,
\end{equation}
where
\begin{enumerate}
\item[$(i)$] $g_1, \, g_2:\mathbb R \longrightarrow \mathbb R$ are nondecreasing and left-continuous everywhere. Moreover, $a\notin N_{g_1}^{-}$ and $b\notin C_{g_1}\cup N_{g_1}^{+}\cup D_{g_1}$;
    \item[$(ii)$] $A:\mathbb{R}\longrightarrow\mathbb{R}$ is positive and continuous;
    \item[$(iii)$] $\Phi:\mathbb{R}\longrightarrow(-c,c)$, $c>0$, is an increasing homeomorphism; 
    \item[$(iv)$] $G:I\times\mathbb{R}^2\longrightarrow\mathcal{P}(\mathbb{R})$ is a multivalued mapping taking non-empty, closed and convex values, and satisfies the following Carathéodory conditions:
    \begin{enumerate}
        \item the mapping $t\to G(t,x,y)$ is $g_2$-measurable for every $(x,y)\in\mathbb{R}^2$, and
        \item the mapping $(x,y)\to G(t,x,y)$ is u.s.c. for 
        $g_2$-almost every $t\in I$.
    \end{enumerate}
\end{enumerate}

\noindent
As in the preceeding section, we define weak and strong solutions for problem \eqref{DP}. First, a {\bf weak} Carathéodory solution of \eqref{DP} is a function $x \in \mathcal{AC}_{g_1}(I)$ such that $x(a)=v_a$, $x(b)=v_b$, and there exists $\xi  \in \mathcal{AC}_{g_2}(I)$ such that
$$\left\{
\begin{array}{l}
\xi(t)=A(x(t)) \Phi (x'_{g_1}(t)) \, \, \mbox{for $g_1$-a.a. $t \in I$, and}\\
\\
\xi'_{g_2}(t) \in G\left(t,x(t),\Phi^{-1}\left(
\dfrac{\xi(t)}{A(x(t))} \right) \right)\, \, \mbox{for $g_2$-a.a. $t \in I$.}
\end{array}
\right.$$
On the other hand, a {\bf strong} Carathéodory solution of \eqref{DP} is a function $x \in \mathcal{BC}^1_{g_1}(I)$ such that $x(a)=v_a$, $x(b)=v_b$, $(A \circ x) \cdot (\Phi \circ x'_{g_1})  \in \mathcal{AC}_{g_2}(I)$, and
$$((A \circ x) \cdot (\Phi \circ x'_{g_1}))'_{g_2}(t)\in G\left(t,x(t),x'_{g_1}(t)  \right) \quad \mbox{for $g_2$-a.a. $t \in I$.}$$

Following the ideas in \cite{FeMaPa}, we will prove the existence of at least one solution for problem \eqref{DP} by combining the fixed-point approach (with Theorem \ref{teo:fun}) and the lower and upper solutions method. In this context of $g$-differential inclusions, lower and upper solutions are defined as follows.
\begin{definition}
    A lower solution for problem \eqref{DP} is a function $\alpha\in\mathcal{BC}_{g_1}^1(I)$ such that $(A\circ \alpha)\cdot(\Phi\circ \alpha'_{g_1})\in\mathcal{AC}_{g_2}(I)$,
    $$
    \alpha(a)\leq v_a,\quad \alpha(b)\leq v_b,
    $$
    and
   \begin{equation}\label{alphag}
   G(t,\alpha(t),\alpha'_{g_1}(t))\cap\left(-\infty,   ((A\circ \alpha) \cdot(\Phi \circ \alpha'_{g_1}))'_{g_2}(t)\right]\neq\emptyset \quad \text{for}\,\, g_2\text{-a.a.}\,\, t \in I.
   \end{equation}

   Analogously, an upper solution for the problem \eqref{DP} is a function $\beta\in\mathcal{BC}_{g_1}^1(I)$ such that $(A\circ \beta)\cdot(\Phi\circ \beta'_{g_1})\in\mathcal{AC}_{g_2}(I)$,
    $$
    \beta(a)\geq v_a,\quad \beta(b)\geq v_b,
    $$
    and
   \begin{equation}\label{betag}
   G(t,\beta(t),\beta'_{g_1}(t))\cap\left[  ((A\circ \beta) \cdot(\Phi \circ \beta'_{g_1}))'_{g_2}(t),+\infty\right)\neq\emptyset \quad \text{for}\,\, g_2\text{-a.a.}\,\, t \in I.
   \end{equation}
\end{definition}

Observe that \eqref{alphag} is equivalent to 
$$G(t,\alpha(t),\alpha'_{g_1}(t)) \le ((A\circ \alpha) \cdot(\Phi \circ \alpha'_{g_1}))'_{g_2}(t)$$
when $G(t,\alpha(t),\alpha'_{g_1}(t))$ is singlevalued. An analogous remark is true for \eqref{betag}.

Given a well-ordered pair of lower and upper solutions $(\alpha,\beta)$, which means that $\alpha\leq\beta$ on $I$, we will denote
$$
\widetilde{I}=[\min_{t\in I}\alpha(t),\max_{t\in I}\beta(t)], \qquad \widetilde{v}=\max_{t\in I}\beta(t)-\min_{t\in I}\alpha(t),
$$
and
$$
m=\min_{x\in\widetilde{I}}A(x)>0,\qquad M=\max_{x\in\widetilde{I}}A(x).
$$
Before moving onto the main theorem of this section, we state and prove the following lemma, adapted from \cite{HuPa2}, Proposition III.2.6, which will be very useful when proving the existence of solution for problem \eqref{DP}.

\begin{lemma}\label{lemaF}
 Let $F:I\times\mathbb{R}^2\rightarrow\mathcal{P}(\mathbb{R})$ a multivalued mapping with non-empty, closed and convex values, verifying
    \begin{enumerate}
                \item[(a)] for every $(x,y)\in\mathbb{R}^2$, the mapping $t\rightarrow F(t,x,y)$ is $g_2$-measurable;
                \item[(b)] for $g_2$-almost every $t\in I$, the mapping $(x,y)\rightarrow F(t,x,y)$ is u.s.c.; and
                \item[(c)] there exists $\gamma\in L^{1}_{g_2}(I)$  such that $\|F(t,x,y)\|\leq\gamma(t)$ for $g_2$-a.a. $t\in I$ and every $(x,y) \in \mathbb R^2$.
    \end{enumerate}
    Define $\widetilde{F}:\mathcal{AC}_{g_1}(I)\rightarrow\mathcal{P}(L^{1}_{g_2}(I))$ as 
    $$
        \widetilde{F}_{x}=\left\{v\in L^{1}_{g_2}(I)\,:\,v(t)\in F(t,x(t),x'_{g_1}(t))\,\,\text{$g_2$-a.e. in $I$}\right\}, \quad x \in \mathcal{AC}_{g_1}(I),
    $$
    where $x'_{g_1}(t)$ is defined as zero whenever $x'_{g_1}(t)$ does not exist. Then, $\widetilde{F}$ takes non-empty, closed and convex values, and it is w-u.s.c.
    \begin{proof}
     Fix $x\in \mathcal{AC}_{g_1}(I)$ and let us see that $\widetilde{F}_x$ is non-empty. Since $x,x'_{g_1}\in L^1_{g_1}(I)$, there exist sequences of simple functions $\{r_n\}_{n\in\mathbb{N}},\{s_n\}_{n\in\mathbb{N}}\subset L^1_{g_1}(I)$ such that $r_n(t)\to x(t)$ and $s_n(t)\to x'_{g_1}(t)$ for $g_1$-a.a. $t\in I$. Hypothesis $(a)$ implies that the mapping $t\to F(t,r_n(t),s_n(t))$ is $g_2$-measurable for every $n\in\mathbb{N}$, so we can apply Kuratowski-Ryll Nardzewski selection theorem, see \cite[Theorem II-2.1]{HuPa1}, to deduce the existence of $v_n\in L^1_{g_2}(I)$ such that $v_n(t)\in F(t,s_n(t),r_n(t))$ for $g_2$-a.a. $t\in I$. From $(c)$, we get $\|v_n\|_{L^1_{g_2}}\leq \|\gamma\|_{L^1_{g_2}}$ for all $n \in \mathbb N$, so we can use the Dunford-Pettis theorem: there exists a function $v\in L^1_{g_2}(I)$ such that $v_n\overset{w}{\to}v$ (passing to a subsequence if necessary). By Mazur's Lemma (see \cite{et}), there exists a sequence $\left\{w_n\right\}_{n \in \mathbb N}$, where each $w_n$ is a convex linear combination of some of the $v_k$'s with $k \ge n$, such that $w_n \to v$ in the $L^1_{g_2}$-norm. Finally, there is a subsequence, which we again denote by $\{ w_n \}_{n \in \mathbb N}$, such that $w_n(t) \to v(t)$ for $g_2$-a.a. $t \in I$. Let $t \in I$ be such that $r_n(t) \to x(t)$, $s_n(t) \to x'_{g_1}(t)$ and $w_n(t) \to v(t)$; we shall prove that $v(t) \in F(t,x(t),x'_{g_1}(t))$, thus proving that $\widetilde{F}_x$ is not empty. Since $F(t,\cdot,\cdot)$ is u.s.c. at $(x(t),x'_{g_1}(t))$, for each $k \in \mathbb N$ there exists $\delta_k>0$ such that
     $$F(t,B((x(t),x'_{g_1}(t));\delta_k))\subset F(t,x(t),x'_{g_1}(t))+(-1/k,1/k),$$
     where $B((x(t),x'_{g_1}(t));\delta_k)$ is the open ball centered at $(x(t),x'_{g_1}(t))$ and radius $\delta_k$. Since $r_n(t) \to x(t)$ and $s_n(t) \to x'_{g_1}(t)$, there exists $n_k$ such that
     $$n \ge n_k \Rightarrow (x_n(t),s_n(t)) \in B((x(t),x'_{g_1}(t));\delta_k),$$
     so for $n \ge n_k$ we have $F(t,x_n(t),s_n(t)) \subset F(t,x(t),x'_{g_1}(t))+(-1/k,1/k)$, which is convex, and therefore
     $$w_{n_k}(t) \in \mbox{co} \, \bigcup_{n=n_k}^{\infty}F(t,x_n(t),s_n(t)) \subset F(t,x(t),x'_{g_1}(t))+(-1/k,1/k).$$
     This implies that $$d(v(t),F(t,x(t),x'_{g_1}(t)))= \lim_{k \to \infty}d(w_{n_k}(t),F(t,x(t),x'_{g_1}(t)))=0,$$
     or, equivalently, $v(t) \in F(t,x(t),x'_{g_1}(t))$ because $F$ assumes closed values.
     
     Proving that  $\widetilde{F}$ takes closed and convex values follows directly from the fact that $F$ does the same.

     \noindent
    Lastly, let us check that $x\to\widetilde{F}_x$ is w-u.s.c. Take $\{x_n\}_{n\in\mathbb{N}}\subset\mathcal{AC}_{g_1}(I)$ such that $x_n\to x$ in $\mathcal{AC}_{g_1}(I)$. We have that $x_n(t)\to x(t)$ and, up to a subsequence, $(x_n)'_{g_1}(t)\to x'_{g_1}(t)$ for $g_1$-almost $t\in I$. Let $v_n\in \widetilde{F}_{x_n}$ for each $n\in\mathbb{N}$ and $v \in L^1_{g_2}(I)$ such that $v_n\overset{w}{\to}{v}$ in $L^1_{g_2}(I)$. Using the same reasoning as before, we get that $v\in\widetilde{F}_x$, so $\tilde{F}$ is w-u.s.c.
        \end{proof}
\end{lemma}

We are now in a position to state and prove the following existence result for the Dirichlet problem.
\begin{theorem}\label{teo:DP}
Let $\alpha$ and $\beta$ be a well-ordered pair of lower and upper solutions for problem \eqref{DP}, and assume that conditions $(i)$-$(iv)$ at the beginning of Section 4 hold. 

Suppose that there exists $\psi\in L_{g_2}(I)$ such that  
\begin{equation}\label{cota2}
    \|G(t,x,y)\|\leq\psi(t)\quad\text{for $g_2$-a.a.}\,\, t\in I\,\,\text{and every} \,\, (x,y)\in\mathbb{R}^2\,\,\text{with}\,\,\alpha(t)\leq x \leq\beta(t),
\end{equation}
and the $g_1$-slope $s^*=(v_b-v_a)/(g_1(b)-g_1(a))$ verifies \eqref{ms}.

Assume also that ${g_2}_{|I} \in {\cal C}_{g_1}(I)$ and

\noindent
either $D_{g_1}=\emptyset$ 

\noindent
or for each fixed $t \in D_{g_1}$ one of the following alternatives hold:
\begin{enumerate}
\item[(a)] $A(\alpha(t)) \ge A(z) \ge A(\beta(t))$ for every $z \in [\alpha(t),\beta(t)]$, and
\begin{enumerate}
    \item[(i)] if $A$ is not constant in $[\alpha(t),\beta(t)]$ and $t \not \in D_{g_2}$, then $\Phi(\alpha'_{g_1}(t)), \Phi(\beta'_{g_1}(t))\ge 0$,
    \item[(ii)] if $t \in D_{g_2}$, then 
    $$G(t,x,y) \subset (-\infty, ((A \circ \alpha)\cdot (\Phi \circ \alpha'_{g_1}))'_{g_2}(t)] \quad \mbox{for $x \in [\alpha(t),\beta(t)]$ and $y \le \alpha'_{g_1}(t)$,}$$
    and
    $$G(t,x,y) \subset [((A \circ \beta)\cdot (\Phi \circ \beta'_{g_1}))'_{g_2}(t),+\infty) \quad \mbox{for $x \in [\alpha(t),\beta(t)]$ and $y \ge \beta'_{g_1}(t)$,}$$
\end{enumerate}
\end{enumerate}
\noindent
or
\begin{enumerate}
\item[(b)] $A(\alpha(t)) \le A(z) \le A(\beta(t))$ for every $z \in [\alpha(t),\beta(t)]$, and
\begin{enumerate}
    \item[(i)] if $A$ is not constant in $[\alpha(t),\beta(t)]$ and $t \not \in D_{g_2}$, then $\Phi(\alpha'_{g_1}(t)), \Phi(\beta'_{g_1}(t))\le 0$,
    \item[(ii)] if $t \in D_{g_2}$, then 
    $$G(t,x,y) \subset (-\infty, ((A \circ \alpha)\cdot (\Phi \circ \alpha'_{g_1}))'_{g_2}(t)] \quad \mbox{for $x \in [\alpha(t),\beta(t)]$ and $y \le \alpha'_{g_1}(t)$,}$$
    and
    $$G(t,x,y) \subset [((A \circ \beta)\cdot (\Phi \circ \beta'_{g_1}))'_{g_2}(t),+\infty) \quad \mbox{for $x \in [\alpha(t),\beta(t)]$ and $y \ge \beta'_{g_1}(t)$,}$$
\end{enumerate}
\end{enumerate}
Then, problem \eqref{DP} has at least one strong solution
$x$ such that $\alpha(t) \le x(t) \le \beta(t)$, $t \in I$. 
    \begin{proof}
        First, we define the truncation operator $T:\mathcal{AC}_{g_1}(I)\rightarrow\mathcal{AC}_{g_1}(I)$ given by
        $$
        T_{x}(t)=\max\{\alpha(t),\min\{x(t),\beta(t) \}\}=
        \left\{
        \begin{array}{lll}
            \alpha(t), & x(t)<\alpha(t),\\
            x(t), & \alpha(t)\leq x(t)\leq\beta(t),\\
            \beta(t), & x(t)>\beta(t).
            
        \end{array}
        \right.
        $$
        Since $\alpha,\beta\in\mathcal{AC}_{g_1}(I)$, Proposition \ref{maxminac} ensures that $T_{x}\in\mathcal{AC}_{g_1}(I)$ for every $x\in\mathcal{AC}_{g_1}(I)$, so $T$ is well-defined. Let us define the auxiliar multivalued mapping $\widetilde{G}:I\times\mathbb{R}^2\longrightarrow\mathcal{P}(\mathbb{R})$ as follows:
        $$
        \widetilde{G}(t,x,y)=
        \left\{
        \begin{array}{lllll}
        G(t,x,y), & \alpha(t)\leq x \leq \beta(t),\\
        C_{y}:=G(t,\alpha(t),y)\cap\left(-\infty,  ((A \circ \alpha) \cdot (\Phi \circ \alpha'_{g_1}))'_{g_2}(t)\right],& x<\alpha(t)\,\,\text{and}\,\, C_{y}\neq\emptyset,\\
        \left\{((A \circ \alpha) \cdot (\Phi \circ \alpha'_{g_1}))'_{g_2}(t)\right\},& x<\alpha(t)\,\,\text{and}\,\, C_{y}=\emptyset,\\ 
        C^{y}:=G(t,\beta(t),y)\cap\left[ ((A \circ \beta) \cdot (\Phi \circ \beta'_{g_1}))'_{g_2}(t),+\infty\right), & x>\beta(t)\,\,\text{and}\,\, C^{y}\neq\emptyset,\\
        \left\{((A \circ \beta) \cdot (\Phi \circ \beta'_{g_1}))'_{g_2}(t)\right\},& x>\beta(t)\,\,\text{and}\,\, C^{y}=\emptyset.
        \end{array}
        \right.
        $$
        Since $G$ takes non-empty, closed, and convex values, so does $\widetilde{G}$. Moreover, it can be easily checked that $\widetilde{G}$ is also Carathéodory.

        \noindent
        Now, note that if $G(t,\alpha(t),y)\cap\left(-\infty,  ((A \circ \alpha) \cdot (\Phi \circ \alpha'_{g_1}))'_{g_2}(t)\right]=\emptyset$, using \eqref{alphag} and \eqref{cota2} we get that
        $$
        \psi(t)\geq\min G(t,\alpha(t),y)\geq((A \circ \alpha) \cdot (\Phi \circ \alpha'_{g_1}))'_{g_2}(t)\geq\min G(t,\alpha(t),\alpha'_{g_1}(t))\geq-\psi(t),
        $$
        for $g_2$-almost all $t\in I$. Analogously, if $G(t,\beta(t),y)\cap\left[ ((A \circ \beta) \cdot (\Phi \circ \beta'_{g_1}))'_{g_2}(t),+\infty\right)=\emptyset$, using \eqref{betag} and \eqref{cota2}, we deduce that for $g_2$-almost $t\in I$
        $$
        -\psi(t)\leq\max G(t,\beta(t),y)\leq ((A \circ \beta) \cdot (\Phi \circ \beta'_{g_1}))'_{g_2}(t)\leq \max G(t,\beta(t),\beta'_{g_1}(t))\leq \psi(t).
        $$
        These inequalities, together with the definition of $\widetilde{G}$ and condition \eqref{cota2}, imply that
        \begin{equation}\label{cota3}
        \|\widetilde{G}(t,x,y)\|\leq\psi(t) \quad\text{for $g_2$-a.a.}\,\, t\in I\,\,\text{and every} \,\, (x,y)\in\mathbb{R}^2.
        \end{equation}

        Now, consider the auxiliary problem
        \begin{equation}\label{plam}
        \left\{
        \begin{array}{ll}
        (\mathcal{A}_x \cdot (\Phi \circ x'_{g_1}))'_{g_2}(t)\in\widetilde{G}(t,x(t),(T_{x})'_{g_1}(t)),\quad t \in I,\\
        x(a)=v_a,\quad x(b)=v_b,
        \end{array}
        \right.
        \end{equation}
     where $\mathcal{A}:\mathcal{AC}_{g_1}(I)\rightarrow\mathcal{BC}_{g_1}(I)$ is given by
        $$
        \mathcal{A}_{x}(t)=A(T_{x}(t)).
        $$
        Note that, since $A$ is continuous, $T$ is bounded and $T_{x}\in\mathcal{AC}_{g_1}(I)$ for every $x\in\mathcal{AC}_{g_1}(I)$, then $A\circ T_{x}\in \mathcal{BC}_{g_1}(I)$ for every $x\in\mathcal{AC}_{g_1}(I)$, so $\mathcal{A}$ is well-defined. Moreover, for every $x\in\mathcal{AC}_{g_1}(I)$ and every $t\in I$ we have that
        $$m\leq \mathcal{A}_{x}(t)\leq M.
        $$
        We will now see that $\mathcal{A}$ is uniformly continuous. Given $\varepsilon>0$, from the uniform continuity of $A$ on $\widetilde{I}$ we deduce the existence of a $\delta>0$ such that $|A(z_1)-A(z_2)|<\varepsilon$ for every $z_1,z_2\in\widetilde{I}$ verifying $|z_1-z_2|<\delta$. Take $x_1,x_2\in\mathcal{AC}_{g_1}(I)$ such that $\|x_1-x_2\|<\delta$. Hence, $\|T_{x_1}-T_{x_2}\|<\delta$ and
        $$
        \|\mathcal{A}_{x_1}-\mathcal{A}_{x_2}\|_{0}=\max|A(T_{x_1}(t))-A(T_{x_2}(t))|<\varepsilon,
        $$
        so $\mathcal{A}$ is uniformly continuous. Lastly, let $\mathcal{S}\subset\mathcal{AC}_{g_1}(I)$ be a bounded set, i.e., there exists $k>0$ such that $\|x\|\leq k$ for every $x\in\mathcal{S}$, and let us check that $\mathcal{A}(\mathcal{S})$ is $g_1$-equicontinuous. For any fixed $\varepsilon>0$, the uniform continuity of $A$ in $[-k,k]$ implies the existence of $\delta>0$ such that $|A(r_1)-A(r_2)|<\varepsilon$ for any $r_1,r_2\in[-k,k]$ such that $|r_1-r_2|<\delta$. Now, taking $\eta=\delta/k$, for any $t,s\in I$ such that $|g_1(t)-g_1(s)|<\eta$ we have
        $$
        \left|T_{x}(t)-T_{x}(s)\right|=\left|\int_{[t,s)}(T_{x})'_{g_1}(\tau)\,dg_1(\tau)\right|\leq\left|\int_{[t,s)}|x'_{g_1}(\tau)|\,d\tau\right|\leq k|g_1(t)-g_1(s)|<\delta,
        $$
        for any $x\in\mathcal{S}$, so
        $$
        |\mathcal{A}_{x}(t)-\mathcal{A}_{x}(s)|=\left|A(T_{x}(t))-A(T_{x}(s))\right|<\varepsilon, \quad \text{for every}\,\,x\in\mathcal{S}.
        $$
        Subsequently, $\mathcal{A}$ maps bounded sets into $g_1$-equicontinuous sets.

        Now, consider the multivalued operator $F:\mathcal{BC}^1_{g_1}(I)\rightarrow\mathcal{P}(L^{1}_{g_2}(I))$ defined by
        $$
        F_x=\left\{
        y \in L^{1}_{g_2}(I) \, : \, y(t) \in\widetilde{G}(t,x(t),(T_{x})'_{g_1}(t)) \,\,\text{$g_2$-a.e. $t\in I$} \right\}.
        $$
         Now, we will use Lemma \ref{lemaF} to deduce that $F$ satisfies the conditions required in Theorem \ref{teo:fun}.
        Notice that $\widetilde{G}$ takes non-empty, closed and convex values, and satisfies $(a)$-$(b)$-$(c)$. Moreover, the proof of Lemma \ref{lemaF} is the same if we replace $x'_{g_1}$ by $(T_{x})'_{g_1}$, so we deduce that $F_{x}$ is non-empty, closed and convex for every $x\in\mathcal{BC}^1_{g_1}(I)$, and $x\to F_x$ is w-u.s.c (because the $\|\cdot\|_{\mathcal{AC}_{g_1}}$ norm is equivalent to the $\|\cdot\|_{1}$ norm in $\mathcal{BC}_{g_1}^1(I)$). Also, from \eqref{cota3} we get 
        $$
        \|F_{x}(t)\|=\sup_{y\in F_x}|y(t)|\leq \psi(t) 
        $$
        for every $x\in\mathcal{BC}^1_{g_1}(I)$ and $g_2$-almost $t\in I$. Hence, $F$ satisfies the hypothesis in Theorem \ref{teo:fun}, so we can apply this result to deduce the existence of a function $x\in\mathcal{AC}_{g_1}(I)$ such that $\mathcal{A}_{x}\cdot(\Phi\circ x'_{g_1})=(A\circ T_{x})\cdot(\Phi\circ x'_{g_1})\in\mathcal{AC}_{g_2}(I)$ and
        $$
        (\mathcal{A}_{x}\cdot(\Phi\circ x'_{g_1}))'_{g_2}(t) \in F_{x}(t)\quad \text{$g_2$-almost all $t\in I$}, \quad x(a)=v_a, \,\, x(b)=v_b.
        $$
        Thus, $x$ is a strong solution to problem \eqref{plam}, and the proof will be complete once we prove that $\alpha(t)\leq x(t)\leq \beta(t)$ for every $t \in I$, which implies that $x$ is a solution of \eqref{DP}. We will only prove that $\alpha(t)\leq x(t)$ for every $t\in I$, since the proof of the other inequality is completely analogous. Suppose, on the contrary, that there exists $t^*\in I$ such that $(x-\alpha)(t^*)<0$. Since $(x-\alpha)(a) \ge 0$, $(x-\alpha)(b) \ge 0$, and $x-\alpha$ is left-continuous everywhere in $(a,b]$ and right-continuous at $t=a$, we deduce the existence of $t_1,t_2 \in [a,b]$, $t_1<t_2$, such that
        \begin{eqnarray}
        \label{signo}
            (x-\alpha)(t)<0 \quad \mbox{for all  $t \in (t_1,t_2)$,} \\
            \label{signot1}
            (x-\alpha)(t_1) \ge 0,\\
          \nonumber
            \mbox{either $t_2=b$ and $(x-\alpha)(b)=0$, or}\\
            \label{signot2} 
            \mbox{$t_2<b$ and $(x-\alpha)(t_2) \le 0 \le (x-\alpha)(t_2^+).$}
        \end{eqnarray} We claim that $x'_{g_1}(t_2) \ge \alpha'_{g_1}(t_2)$. Indeed, if $t_2<b$ and $t_2 \in D_{g_1}$, then condition (\ref{signot2}) and the definition of $g_1$-derivative at discontinuity points immediately yield $(x-\alpha)'_{g_1}(t_2) \ge 0$. Assume, on the other hand, that $t_2=b$ or $t_2<b$ and $t_2 \not \in D_{g_1}$. In both cases, $x-\alpha$ is continuous at $t_2$, so $(x-\alpha)(t_2)=0$. Moreover, we have $g_1(t)<g_1(t_2)$ for $t<t_2$. To prove it, observe that if for some $\varepsilon >0$ we had $g_1(t)=g_1(t_2)$ for $t \in (t_2-\varepsilon,t_2]$, then $x-\alpha$ would be constant in $(t_2-\varepsilon,t_2]$ (because $x-\alpha \in \mathcal{BC}_{g_1}(I)$) but this would contradict \eqref{signo} and $(x-\alpha)(t_2)=0$. Finally, there exists $(x-\alpha)'_{g_1}(t_2)$ because $x-\alpha \in \mathcal{BC}^1_{g_1}(I)$, hence
        $$(x-\alpha)'_{g_1}(t_2)=\lim_{t \to t_2^-}\dfrac{(x-\alpha)(t)}{g_1(t)-g_1(t_2)} \ge 0,$$
        and the claim is proven.

        \bigbreak
        
        Now, since $x$ is a solution of \eqref{plam} and $x-\alpha< 0$ in $(t_1,t_2)$, we have $T_{x}(t)=\alpha(t)$ in $(t_1,t_2)$ and  
        $$
        ((A \circ \alpha) \cdot (\Phi\circ x'_{g_1}))'_{g_2}(t) \in G(t,\alpha(t),\alpha'_{g_1}(t))\cap\left(-\infty,  ((A \circ \alpha) \cdot (\Phi\circ \alpha'_{g_1}))'_{g_2}(t)\right],
        $$
        for $g_2$-almost all $t\in (t_1,t_2)$. Hence, 
        $$
        ((A \circ \alpha) \cdot (\Phi\circ x'_{g_1}))'_{g_2}(t)-((A \circ \alpha) \cdot (\Phi\circ \alpha'_{g_1}))'_{g_2}(t)\leq 0 \quad \text{for $g_2$-a.a. $t\in (t_1,t_2)$.}$$ 
        For every $t \in (t_1,t_2)$ and every $s \in (t,t_2)$, the function $\Theta(r)=A(\alpha(r))(\Phi(x'_{g_1}(r))-\Phi(\alpha'_{g_1}(r)))$ is $g_2$-absolutely continuous in $[t,s]$, and therefore we have
        $$A(\alpha(s))(\Phi(x'_{g_1}(s))-\Phi(\alpha'_{g_1}(s)))-A(\alpha(t))(\Phi(x'_{g_1}(t))-\Phi(\alpha'_{g_1}(t)))=\int_{[t,s)} \Theta'_{g_2}(r) \, dg_2(r) \le 0,$$
        thus
        \begin{equation}
        \label{otdes}
        A(\alpha(t))(\Phi(x'_{g_1}(t))-\Phi(\alpha'_{g_1}(t)))\ge A(\alpha(s))(\Phi(x'_{g_1}(s))-\Phi(\alpha'_{g_1}(s))) \quad \mbox{whenever $t_1<t<s<t_2$.}
        \end{equation}
        Since $\Theta$ is left-continuous at $t_2$, we deduce that for every $t\in(t_1,t_2)$ we have
        \begin{align*}
        A(\alpha(t))(\Phi(x'_{g_1}(t))-\Phi(\alpha'_{g_1}(t)))&\geq \lim_{s \to t_2^-}A(\alpha(s))(\Phi(x'_{g_1}(s))-\Phi(\alpha'_{g_1}(s)))\\
        & =A(\alpha(t_2))(\Phi(x'_{g_1}(t_2))-\Phi(\alpha'_{g_1}(t_2))) \ge 0,
        \end{align*}
        so $\Phi(x'_{g_1}(t))-\Phi(\alpha'_{g_1}(t))\geq0$ on $(t_1,t_2)$ by the positiveness of $D$, and \begin{equation}
            \label{primapos}
            (x-\alpha)'_{g_1}(t)\geq 0 \quad \mbox{for all $t \in (t_1,t_2)$,} 
        \end{equation}by the monotonicity of $\Phi$. We shall prove that 
       \begin{equation}
           \label{ultcond}
           \mbox{if $t_1 \in D_{g_1}$, then $(x-\alpha)'_{g_1}(t_1)=0$,}
       \end{equation} 
  and then, for both cases $t_1 \in D_{g_1}$ or $t_1 \not \in D_{g_1}$, for any $t\in(t_1,t_2)$ we shall have
        $$
        (x-\alpha)(t)=(x-\alpha)(t_1)+\int_{[t_1,t)}(x-\alpha)'_{g_1}(s)\,dg_1(s)\ge \int_{(t_1,t)}(x-\alpha)'_{g_1}(s)\,dg_1(s)\geq0,
        $$
a contradiction with (\ref{signo}), thus completing the proof that $\alpha \le x$ in $I$.  

\bigbreak
        To begin with the proof of (\ref{ultcond}), observe that inequalities (\ref{signo}) and (\ref{signot1}) ensure that 
        \begin{equation}
\label{xmalfaprim}
x'_{g_1}(t_1) \le \alpha'_{g_1}(t_1).
        \end{equation}
        Also observe that (\ref{primapos}) implies that
        \begin{equation}
\label{xmalfa}
x'_{g_1}(t_1^+) \ge \alpha'_{g_1}(t_1^+).
        \end{equation}
        
        Next we assume that conditions $(a)$ hold. The proof under conditions $(b)$ is similar and we omit it. We consider two cases separately.

        \bigbreak
        
        \noindent
        {\it Case 1. $t_1 \in D_{g_1} \setminus D_{g_2}$.} By definition of solution of \eqref{plam}, we have that $\xi=(A \circ T_x) \cdot (\Phi \circ x'_{g_1}) \in \mathcal{AC}_{g_2}(I),$ and $\xi$ is continuous from the right at $t_1$ because $t_1 \not \in D_{g_2}$. Since $T_x(t)=\alpha(t)$ for $t \in (t_1,t_2)$, we have
        $$x'_{g_1}(t_1^+)=\Phi^{-1} \left(\dfrac{\xi(t_1)}{A(\alpha(t_1^+))} \right).$$
Analogously, $\eta=(A \circ \alpha) \cdot (\Phi \circ \alpha'_{g_1}) \in \mathcal{AC}_{g_2}(I)$ and  
$$\alpha'_{g_1}(t_1^+)=\Phi^{-1} \left(\dfrac{\eta(t_1)}{A(\alpha(t_1^+))} \right).$$
Therefore, the inequality (\ref{xmalfa}) yields $\xi(t_1) \ge \eta(t_1)$, thus
$$
A(T_x(t_1))\Phi(x'_{g_1}(t_1)) \ge A(\alpha(t_1))\Phi(\alpha'_{g_1}(t_1)).
$$
The condition $(a)(i)$ ensures that $A$ is constant or $A(\alpha(t_1))\ge A(T_x(t_1))$ and $\Phi(\alpha'_{g_1}(t_1)) \ge 0$, so in both cases we have
$$A(T_x(t_1))\Phi(x'_{g_1}(t_1))  \ge A(T_x(t_1))\Phi(\alpha'_{g_1}(t_1)).$$Hence, obviously, $x'_{g_1}(t_1) \ge \alpha'_{g_1}(t_1)$, which combined with \eqref{xmalfaprim} produces $(x-\alpha)'_{g_1}(t_1) =0$ as desired.  
        \bigbreak
        
        \noindent
        {\it Case 2. $t_1 \in D_{g_1} \cap D_{g_2}$.} Since $t_1\in D_{g_2}$, it has positive $g_2$-measure, so we deduce from \eqref{alphag} that 
        \begin{equation*}\label{ay}
           G(t,T_x(t_1),(T_x)'_{g_1}(t_1))\cap\left[((A\circ \beta) \cdot(\Phi \circ \beta'_{g_1}))'_{g_2}(t_1),+\infty\right)\neq\emptyset \quad \text{if}\quad x(t_1)>\beta(t_1),
        \end{equation*}
        \noindent
        since $(T_x)'_{g_1}(t)=\alpha'_{g_1}(t)$ whenever $x(t)<\alpha(t)$. Analogously, from \eqref{betag} we deduce that
        \begin{equation*}\label{by}
        G(t,T_x(t_1),(T_x)'_{g_1}(t_1))\cap\left(-\infty,   ((A\circ \alpha) \cdot(\Phi \circ \alpha'_{g_1}))'_{g_2}(t_1)\right]\neq\emptyset \quad \text{if}\quad x(t_1)>\beta(t_1),
        \end{equation*}
        
        since $(T_x)'_{g_1}(t)=\beta'_{g_1}(t)$ whenever $x(t)>\beta(t)$. Hence, we have that
        \begin{equation*}
            \widetilde{G}(t_1,x(t_1),(T_x)'_{g_1}(t_1))=\left\{
            \begin{array}{lll}
             G(t_1,T_x(t_1),(T_x)'_{g_1}(t_1)),& \alpha(t)\leq x(t_1) \leq \beta(t_1),\\
             G(t_1,T_x(t_1),(T_x)'_{g_1}(t_1))\cap\left(-\infty,   ((A\circ \alpha) \cdot(\Phi \circ \alpha'_{g_1}))'_{g_2}(t_1)\right], & x(t_1)<\alpha(t_1),\\
             G(t_1,T_x(t_1),(T_x)'_{g_1}(t_1))\cap\left[((A\circ \beta) \cdot(\Phi \circ \beta'_{g_1}))'_{g_2}(t_1),+\infty\right),& x(t_1)>\beta(t_1).\\
            \end{array}
            \right.
        \end{equation*}
        Now, since $x$ is a solution of \eqref{plam}, we get that
        $$((A \circ T_x) \cdot (\Phi\circ x'_{g_1}))'_{g_2}(t_1) \in \widetilde{G}(t_1,x(t_1),(T_x)'_{g_1}(t_1))\subset G(t_1,T_x(t_1),(T_x)'_{g_1}(t_1)).$$
        We deduce from \eqref{signo} and \eqref{signot1} that 
        \begin{align*}
            (T_x)'_{g_1}(t_1)&=\dfrac{T_x(t_1^+)-T_x(t_1)}{g_1(t_1^+)-g_1(t_1)}=\dfrac{\alpha(t_1^+)-T_x(t_1)}{g_1(t_1^+)-g_1(t_1)} \\
            & \le \dfrac{\alpha(t_1^+)-\alpha(t_1)}{g_1(t_1^+)-g_1(t_1)}=\alpha'_{g_1}(t_1),
        \end{align*}
        and then condition $(a)(ii)$ ensures that
        $$
        ((A \circ T_x) \cdot (\Phi\circ x'_{g_1}))'_{g_2}(t_1)  \le ((A \circ \alpha) \cdot (\Phi\circ \alpha'_{g_1}))'_{g_2}(t_1).
      $$
       By definition of $g_2$-derivative,
       $$A(T_x(t_1^+))\Phi(x'_{g_1}(t_1^+))-A(T_x(t_1))\Phi(x'_{g_1}(t_1)) \le A(\alpha(t_1^+))\Phi(\alpha'_{g_1}(t_1^+))-A(\alpha(t_1))\Phi(\alpha'_{g_1}(t_1)).$$
    Notice that (\ref{signo}) guaranties that $T_x(t_1^+)=\alpha(t_1^+)$, so \eqref{xmalfaprim} and conditions $(a)$ yield         \begin{align}
         \nonumber
         A(\alpha(t_1^+))\left(\Phi(x'_{g_1}(t_1^+))- \Phi(\alpha'_{g_1}(t_1^+)) \right)&\le  A(T_x(t_1))\Phi(x'_{g_1}(t_1))- A(\alpha(t_1))\Phi(\alpha'_{g_1}(t_1)) \\
         \label{Phiprima}
         & \le \Phi(\alpha'_{g_1}(t_1)) \left(A(T_x(t_1))-A(\alpha(t_1)) \right) \le 0.
         \end{align} 
         On the other hand, (\ref{xmalfa}) implies that
         $$0 \le A(\alpha(t_1^+))\left(\Phi(x'_{g_1}(t_1^+))- \Phi(\alpha'_{g_1}(t_1^+)) \right),$$
         which, combined with (\ref{Phiprima}), yield
         \begin{equation}
\label{eqno1}
A(T_x(t_1))\Phi(x'_{g_1}(t_1))- A(\alpha(t_1))\Phi(\alpha'_{g_1}(t_1))=0
\end{equation}
\noindent
and
\begin{equation}
         \label{eqno2}
        \Phi(\alpha'_{g_1}(t_1)) \left(A(T_x(t_1))-A(\alpha(t_1)) \right)=0.
         \end{equation}
         Let us prove that \eqref{eqno1} and \eqref{eqno2} imply that $(x-\alpha)'_{g_1}(t_1)=0$. We deduce from \eqref{eqno2} that $\Phi(\alpha'_{g_1}(t_1))=0$ or $A(T_x(t_1))-A(\alpha(t_1))=0$. First, if $\Phi(\alpha'_{g_1}(t_1))=0$, then we deduce from \eqref{eqno1} that $\Phi(x'_{g_1}(t_1))=0$, and therefore $x'_{g_1}(t_1)=\alpha'_{g_1}(t_1)$. Second, if $A(T_x(t_1))=A(\alpha(t_1))$, then we get $(x-\alpha)'_{g_1}(t_1)=0$ from \eqref{eqno1}.  
\end{proof}
\end{theorem}

In the following example, we will use this theorem to prove the existence of solution for a particular problem.

\begin{example}
    Consider the following differential inclusion
    \begin{equation}\label{ex1}
       \left( (2+\sin(x))\cdot\arctan(x'_{g_1}-1)\right)'_{g_2}(t)\in G(t, x(t), x'_{g_1}(t)),\quad t\in[0,1],\quad x(0)=0,\,\,x(1)=1,
    \end{equation}
    where $g_1(t)=t+0.5H(t-0.5)$, $g_2(t)=t+H(t-0.5)$ and
    \begin{equation*}
        G(t, x, y) = \begin{cases}
        \{-0.1\mu x^2\}, &  y > 0, \\
        [-0.1\mu x^2, 0.1\mu x^2], &  y = 0, \\
        \{0.1\mu x^2\}, &  y < 0,
    \end{cases}
    \end{equation*}
    where $\mu>0$ is a given parameter.
    We will use see Theorem \ref{teo:DP} to prove that, for sufficiently small $\mu$,  problem \eqref{ex1} has at least one solution.

    We will first verify that hypotheses $(i)$-$(iv)$ are satisfied. It is obvious that $A(x)=2+\sin(x)$ is positive ($1\leq A(x)\leq3$) and continuous, and $\Phi(s)=\arctan(s-1)$ is an increasing homeomorphism mapping $\mathbb{R}$ onto $(-\pi/2,\pi/2)$. It is also immediate that $G(t,x,y)$ takes non-empty, closed, and convex values. Moreover, for each fixed $(x,y)\in\mathbb{R}^2$, the mapping $t\to G(t,x,y)$ is constant, hence measurable, and for any fixed $t\in I$ the mapping $(x,y)\to G(t,x,y)$ is u.s.c., since its graph is closed and it assumes compact values. To see the graph is closed, take a sequence $\{(x_n,y_n,z_n)\}_n\in\mathrm{Gr}(G(t,\cdot,\cdot))$ such that $(x_n,y_n,z_n)\to(x,y,z)$ and let us see that $(x,y,z)\in\mathrm{Gr}(G(t,\cdot,\cdot))$, i.e., that $z\in G(t,x,y)$. Assume $y=0$ and $y_n\neq 0$ for infinitely many $n\in \mathbb{N}$ (in any other case, the proof is trivial). Assume also, without loss of generality, that $y_n<0$ for every $n\in \mathbb{N}$. Hence, $z_n=0.1\mu x_n^2$ for every $n\in\mathbb{N}$ and $z=\lim_{n\to\infty}0.1\mu x_n^2=0.1\mu x^2\in[-0.1\mu x^2,0.1\mu x^2]=G(t,x,0)=G(t,x,y)$.

    Let us now check that, for a sufficiently small $\mu>0$, the functions 
    $$ \alpha(t) =- \dfrac{2}{3}t - \dfrac{1}{3}H(t-0.5) \quad \text{and} \quad \beta(t) = \dfrac{2}{3}t + \dfrac{1}{3}H(t-0.5) $$
    form a well-ordered pair of lower and upper solutions for the problem \eqref{ex1}. First, computing their $g_1$-derivative, we get $\alpha'_{g_1}(t) = -2/3$ and $\beta'_{g_1}(t) = 2/3$ for all $t \in I$. Since constant functions are $g_1$-continuous, $\alpha, \beta \in \mathcal{BC}_{g_1}^1(I)$. To check that $(A\circ\alpha)\cdot(\Phi\circ\alpha'_{g_1})\in\mathcal{AC}_{g_2}(I)$, observe first that $\Phi(\alpha'_{g_1}(t))=\arctan(-5/3)$ is constant for every $t\in I$, so it suffices to check that  $(A\circ\alpha)\in\mathcal{AC}_{g_2}(I)$, which follows from the fact that $A$ is a Lipschitz-continuous function and $\alpha\in\mathcal{AC}_{g_2}(I)$ (because $\alpha\in\mathcal{AC}_{g_1}(I)$ and $g_2$ is $g_1$-continuous). In a completely analogous way, one proves that $(A\circ\beta)(\Phi\circ\beta'_{g_1})\in\mathcal{AC}_{g_2}(I)$. Moreover, some easy computations show that
    \begin{equation*}
    ((A\circ\alpha)\cdot(\Phi\circ\alpha'_{g_1}))'_{g_2}(t) = 
    \begin{cases}
        -2/3 \Phi(-2/3) \cos(\alpha(t)), & t \neq 0.5, \\
        \Phi(-2/3) [\sin(-2/3) - \sin(-1/3)], & t = 0.5,
        \end{cases}
    \end{equation*}
    and
    \begin{equation*}
        ((A\circ\beta)(\Phi\circ\beta'_{g_1}))'_{g_2}(t) = 
            \begin{cases}
            2/3 \Phi(2/3) \cos(\beta(t)), & t \neq 0.5, \\
            \Phi(2/3) [\sin(2/3) - \sin(1/3)], & t = 0.5.
        \end{cases}
    \end{equation*}
    Since $\alpha(t)\in[-1,0]$ and $\beta(t)\in[0,1]$, we deduce that
    $$
    ((A\circ\alpha)\cdot(\Phi\circ\alpha'_{g_1}))'_{g_2}(t)\geq0.3711\quad\text{and}\quad((A\circ\beta)(\Phi\circ\beta'_{g_1}))'_{g_2}(t)<-0.1159
    ,
    $$
    for every $t\in I$, $t\neq0.5$.
    \noindent
    Furthermore, since $\alpha'_{g_1}(t)=-0.5$, we have  $G(t,\alpha(t),\alpha'_{g_1})=\{0.1\mu(\alpha(t))^2\}$ for every $t\in I$, so
    $$
    G(t,\alpha(t),\alpha'_{g_1})\cap(-\infty,((A\circ\alpha)\cdot(\Phi\circ\alpha'_{g_1}))'_{g_2}(t)]\neq\emptyset,
    $$
    provided that 
    $$
    \mu\leq \dfrac{0.3711}{0.1\cdot1^2}=3.711\quad \text{and}\quad\mu\leq \dfrac{\Phi(-2/3) [\sin(-2/3) - \sin(-1/3)]}{0.1\cdot\alpha(0.5)^2}=27.
    $$
    Analogously, we get
    $$
    G(t,\beta(t),\beta'_{g_1})\cap[(A\circ\beta)\cdot(\Phi\circ\beta'_{g_1}))'_{g_2}(t),+\infty)\neq\emptyset
    $$
    for every $t\in I$ whenever
    $$
    \mu\leq \dfrac{0.1159}{0.1\cdot1^2}=1.159\quad \text{and}\quad\mu\leq -\dfrac{\Phi(2/3) [\sin(2/3) - \sin(1/3)]}{0.1\cdot\beta(0.5)^2}=8.433.
    $$
    Hence, we conclude that $\alpha$ and $\beta$ form a well ordered pair of lower and upper solutions for the problem \eqref{ex1} for every $\mu\in(0,1.159]$, since $\alpha(0)=0=\beta(0)$, $\alpha(1)=-1\leq1$, $\beta(1)=1$ and $\alpha\leq\beta$ in $I$.

    Now, observe that 
    $$
    \|G(t,x,y)\|\leq0.1\mu\quad\text{for every}\,\, t\in I\,\,\text{and every} \,\, (x,y)\in\mathbb{R}^2\,\,\text{with}\,\,\alpha(t)\leq x \leq\beta(t),
    $$
    and, since the function $\psi(t)=0.1\mu$ trivially belongs to $\mathcal{L}_{g_2}^1(I)$, condition \eqref{cota2} is satisfied. Moreover, $\|\psi\|_{\mathcal{L}_{g_2}^1}=\int_{[a,b)}0.1\mu\,dg_2(s)=0.2\mu$, so
    $$
    3\left|\Phi\left(\dfrac{1-0}{g_1(1)-g_1(0)}\right)\right|=3|\arctan(-1/3)|=0.9653<1.1072=\dfrac{\pi}{2}-0.4\cdot1.159\leq\dfrac{\pi}{2}-2\|\psi\|_{\mathcal{L}_{g_2}^1},
    $$
    and condition \eqref{ms} is also fulfilled.

    Finally, we will see that condition (b) in Theorem \ref{teo:DP} also holds. 
    
    \noindent 
    On the one hand, since $\alpha(t),\beta(t)\in[-1,1]$ for every $t\in I$ and $A(x)=2+\sin(x)$ is increasing in such interval, it follows that $A(\alpha(t)) \le A(z) \le A(\beta(t))$ for every $z \in [\alpha(t),\beta(t)]$. Moreover, $\Phi(\alpha'_{g_1}(t))=\arctan(-5/3)\leq 0$ and $\Phi(\beta'_{g_1}(t))=\arctan(-1/3)\leq 0$ for every $t\in I$, so part $(i)$ of condition (b) is satisfied.

    \noindent
    On the other hand, in order to check part $(ii)$ observe that $D_{g_2}=\{0.5\}$. For $x\in[\alpha(0.5),\beta(0.5)]=[-1/3,1/3]$ and $y\leq\alpha'_{g_1}(0.5)=-2/3<0$, we have that $G(0.5,x,y)=\{0.1\mu x^2\}$. Since $((A\circ\alpha)\cdot(\Phi\circ\alpha'_{g_1}))'_{g_2}(0.5)=0.3$ and $0.1\mu x^2\leq 0.0387$ for every $x\in[-1/3,1/3]$ and $\mu\in(0,1.159]$, we get that
    $$
    G(0.5,x,y) \subset (-\infty, ((A \circ \alpha)\cdot (\Phi \circ \alpha'_{g_1}))'_{g_2}(0.5)] \quad \text{for $x \in [\alpha(0.5),\beta(0.5)]$ and $y \le \alpha'_{g_1}(0.5)$}.
    $$
    Analogously, for $x\in[-1/3,1/3]$ and $y\geq\beta'_{g_1}(0.5)=2/3>0$, we have that $G(0.5,x,y)=\{-0.1\mu x^2\}$. Since $((A\circ\beta)\cdot(\Phi\circ\beta'_{g_1}))'_{g_2}(0.5)=-0.0937$ and $-0.1\mu x^2\geq-0.0128$ for every $x\in[-1/3,1/3]$ and $\mu\in(0,1.159]$, we get that
    $$
    G(0.5,x,y) \subset [((A \circ \beta)\cdot (\Phi \circ \beta'_{g_1}))'_{g_2}(0.5),+\infty) \quad \text{for $x \in [\alpha(0.5),\beta(0.5)]$ and $y \ge \beta'_{g_1}(0.5)$}.
    $$
    Hence, part $(ii)$ of condition (b) is satisfied, so we can apply Theorem \ref{teo:DP} to deduce the existence of a least one solution $x\in\mathcal{BC}^1_{g_1}(I)$ for problem \eqref{ex1}, verifying $\alpha(t)\leq x(t)\leq \beta(t)$ for every $t\in I$.
\end{example}

 Although the uniqueness of the solutions for problem \eqref{DP} is not guaranteed, we can establish a result regarding the continuous dependence of the solutions on the boundary data.
 
 \begin{lemma}
\label{lem:DC}
    In the conditions of Theorem \ref{teo:DP}, assume that $\alpha(a)\le v_n^a \le \beta(a)$ and $\alpha(b) \le v_n^b  \le \beta(b)$ for all $n \in \mathbb N$, and $v_n^a\to v_0^a$, $v_n^b \to v_0^b$. Denote
    $$
    s_n=\dfrac{v_n^b-v_n^a}{g_1(b)-g_1(a)},\quad s_0=\dfrac{v_0^b-v_0^a}{g_1(b)-g_1(a)},
    $$
    and assume condition \eqref{ms} is satisfied for $s^*=s_n$, for every $n\geq 0$. 

    Let $\{x_n\}_{n\in\mathbb{N}}$ be a sequence in $\mathcal{BC}^1_{g_1}(I)$ such that $x_n$ is a solution to problem \eqref{DP} with boundary conditions $v_a=v_n^a$ and $v_b=v_n^b$, verifying $\alpha(t)\leq x_n(t)\leq\beta(t)$ for every $t\in I$. Then, there exists a subsequence of $\{x_n\}_{n\in\mathbb{N}}$ converging uniformly to $x$, a solution of problem \eqref{DP} with boundary conditions $v_a=v_0^a$ and $v_b=v_0^b$. Moreover, the subsequence of $g_1$-derivatives converges pointwise to $x'_{g_1}$ in I.
\end{lemma}
    \begin{proof}
By definition of solution, for each $n \in \mathbb N$ we have
\begin{equation}
    \label{equisene}
x_n(t)=v_n^a+\int_{[a,t)}{\Phi^{-1} \left( \dfrac{\beta_n+\int_{[a,s)}{f_n(r) \, dg_2(r)}}{A(x_n(s))} \right) \, dg_1(s)}  \quad \mbox{for all $t \in I$,}
\end{equation}
where $f_n \in L_{g_2}^1(I)$, $f_n(t) \in G(t,x_n(t),(x_n)'_{g_1}(t))$ for $g_2$-a.a. $t \in I$, and $\beta_n$ is the unique solution of
\begin{equation}
\label{equisnbeta}
\int_{[a,b)}{\Phi^{-1} \left( \dfrac{\beta_n+\int_{[a,s)}{f_n(r) \, dg_2(r)}}{ A(x_n(s))} \right) \, dg_1(s)}=v_n^b-v_n^a.
\end{equation} 
Since ${g_2}_{|I} \in \mathcal{C}_{g_1}(I)$, the integrand in \eqref{equisene} belongs to $\mathcal{BC}_{g_2}(I) \subset \mathcal{BC}_{g_1}(I)$, and then the second part in Proposition \ref{tffemato} guaranties that
\begin{equation}
        \label{xnprim}
        (x_n)'_{g_1}(t)=\Phi^{-1}\left(\dfrac{\beta_n+\int_{[a,t)}{f_n(s) \, dg_2(s)}}{A(x_n(t))} \right)  \quad \mbox{for all $t \in I$.}
        \end{equation}
Similar ideas to those in the proof of Theorem \ref{teo:fun}, Claim 2, ensure that
\begin{equation}
    \label{eqbetan}
|\beta_n| \le M |\Phi(s_n)|+\|\psi\|_{L^1_{g_2}},
\end{equation}
so for all $t \in I$ we have
\begin{align*}
    \nonumber
    |\Phi((x_n)'_{g_1}(t))|&=\left|\dfrac{\beta_n+\int_{[a,t)}{f_n(s) \, dg_2(s)}}{A(x_n(t))}\right| \le \dfrac{M|\Phi(s_n)|+2\|\psi\|_{L^1_{g_2}}}{m} \\
    &
    \le c-\dfrac{\|\psi\|_{L^1_{g_2}}}{m}\equiv\tilde{c}, \quad \mbox{thanks to (\ref{ms}) for $s^*=s_n$.}
    \end{align*}
    In particular, for all $n \in \mathbb N$ we have
    \begin{equation}
        \label{xprimcota}
        (x_n)'_{g_1}(t)  \in \Phi^{-1}([-\tilde{c},\tilde{c}])\equiv K \quad \mbox{for all $t \in I$,}
        \end{equation}
        where $K$ is a compact subset of the reals. Also, $x_n(a)=v_n^a \in [\alpha(a),\beta(a)]$, $n \in \mathbb N$, and therefore $\{x_n\}
_{n \in \mathbb N}$ is a relatively compact sequence in $\mathcal{BC}_{g_1}(I)$ by Theorem \ref{relcomp}. Therefore, we may assume that $x_n \to x \in \mathcal{BC}_{g_1}(I)
$ uniformly in $I$, so, in particular, $x(a)=v^a$ and $x(b)=v^b$. We shall prove that $x$ solves the differential equation in \eqref{DP} and $(x_n)'_{g_1}(t) \to x'_{g_1}(t)$ for all $t \in I$.

Condition \eqref{cota2} implies that 
$\{f_n\}_{n \in \mathbb N}$ is a bounded sequence in $L^1_{g_2}(I)$, so we may assume that 
$f_n \overset{w}{\to} f \in L^1_{g_2}(I)$. We may also assume that $\beta_n \to \beta_{\infty}$, so for each fixed $t \in I$ we have
\begin{equation}
    \label{limxprim}
(x_n)'_{g_1}(t)=\Phi^{-1}\left(\dfrac{\beta_n+\int_{[a,t)}{f_n(s) \, dg_2(s)}}{A(x_n(t))} \right) \to \Phi^{-1}\left(\dfrac{\beta_{\infty}+\int_{[a,t)}{f(s) \, dg_2(s)}}{A(x(t))} \right).
\end{equation}
Moreover, \eqref{xprimcota} ensures that we can apply the dominated convergence theorem in \eqref{equisene} to deduce
$$x(t)=v_0^a+\int_{[a,t)}{\Phi^{-1} \left( \dfrac{\beta_{\infty}+\int_{[a,s)}{f(r) \, dg_2(r)}}{A(x(s))} \right) \, dg_1(s)}  \quad \mbox{for all $t \in I$,}
$$
which, in particular, implies that $x \in \mathcal{BC}_{g_1}^1(I)$ and
\begin{equation}
    \label{xprimm}
x'_{g_1}(t)=\Phi^{-1}\left(\dfrac{\beta_{\infty}+\int_{[a,t)}{f(s) \, dg_2(s)}}{A(x(t))} \right)  \quad \mbox{for all $t \in I$.}
\end{equation}
By \eqref{limxprim}, $(x_n)'_{g_1}(t) \to x'_{g_1}(t)$ for all $t \in I$, and by \eqref{xprimm}, $(A \circ x)(\Phi \circ x'_{g_1}) \in \mathcal{AC}_{g_2}(I)$ and $((A \circ x)(\Phi \circ x'_{g_1}))'_{g_2}(t)=f(t)$ for $g_2$-a.a. $t \in I$.

Finally, to prove that $f(t) \in G(t,x(t),x'_{g_1}(t))$ for $g_2$-a.a. $t \in I$ it suffices to repeat the final part of the proof of Lemma \ref{lemaF}. 
\end{proof}

\begin{remark}
In the conditions of Lemma \ref{lem:DC} it can be proven that the subsequence of $g_1$-derivatives also converges uniformly in $I$ to the $g_1$-derivative of the solution.
\end{remark}

\section{Periodic and nonlinear boundary conditions}\label{sec5}
In this section, we will prove the existence of solution for the inclusion \eqref{eqg} with different types of general non-linear boundary conditions, covering the particular cases of the periodic, Neumann and Robin boundary conditions.

First, we will consider the following problem
\begin{equation}\label{NLP}
    \left\{
    \begin{array}{lll}
   ((A \circ x)\cdot (\Phi \circ x'_{g_1}))'_{g_2}(t)\in G(t,x(t),x'_{g_1}(t)) \quad \mbox{for $g_2$-a.a. $t \in I=[a,b],$}
   \\ k(x(a),x(b),x'_{g_1}(a),x'_{g_1}(b))=0,\\
   x(b)=h(x(a)),
   \end{array}
   \right.
\end{equation}
where $h:\mathbb{R}\rightarrow\mathbb{R}$ is continuous and increasing, and $k:\mathbb{R}^4\rightarrow\mathbb{R}$ is continuous and verifies that 
\begin{equation}\label{condK}
z\to k(x,y,z,w)\,\,\, \text {is increasing}\quad \text{and} \quad w\to k(x,y,z,w) \, \, \,\text{is decreasing}.  
\end{equation}

\noindent
Our existence result for this kind of problems relies again on the existence of a well-ordered pair of lower and upper solutions, which in this context are defined as follows.

\begin{definition}
    A lower solution for the problem \eqref{NLP} is a function $\alpha\in\mathcal{BC}^1_{g_1}(I)$ such that 
    $$
    k(\alpha(a),\alpha(b),\alpha'_{g_1}(a),\alpha'_{g_1}(b))\geq 0,\quad \alpha(b)=h(\alpha(a)),
    $$
 and $(A\circ x)\cdot(\Phi\circ \alpha'_{g_1})\in\mathcal{AC}_{g_2}(I)$ and
   \begin{equation*}
   G(t,\alpha(t),\alpha'_{g_1}(t))\cap\left(-\infty,  ((A \circ \alpha)\cdot (\Phi \circ \alpha'_{g_1}))'_{g_2}(t)\right]\neq\emptyset \quad \text{for}\,\, g_2\text{-a.a.}\,\, t \in [a,b].
   \end{equation*}

   Analogously, an upper solution for the problem \eqref{NLP} is a function $\beta\in\mathcal{BC}^1_{g_1}(I)$ such that 
    $$
    k(\beta(a),\beta(b),\beta'_{g_1}(a),\beta'_{g_1}(b))\leq 0,\quad \beta(b)=h(\beta(a)),
    $$
    and $(A\circ \beta)\cdot(\Phi\circ \beta'_{g_1})\in\mathcal{AC}_{g_2}(I)$ and
   \begin{equation*}
   G(t,\beta(t),\beta'_{g_1}(t))\cap\left[  ((A \circ \beta)\cdot (\Phi \circ \beta'_{g_1}))'_{g_2}(t),+\infty\right)\neq\emptyset, \quad \text{for}\,\, g_2\text{-almost}\,\, t \in I.
   \end{equation*}
\end{definition}

\noindent
Based on the results for the Dirichlet problem proved on the previous section, we get the following existence result for problem \eqref{NLP}.

\begin{theorem}\label{teo:NLP}
 In the conditions of Theorem \ref{teo:DP} with $\alpha$ and $\beta$   a well-ordered pair of lower and upper solutions for problem \eqref{NLP}, assume that condition \eqref{condK} holds and assume also that the estimate \eqref{ms} is satisfied for $s^*=(h(\eta)-\eta)/(g_1(b)-g_1(a))$, for every $\eta\in[\alpha(a),\beta(a)]$.

 \noindent
 Then, problem \eqref{NLP} has at least one solution $x(t)$ such that $\alpha(t)\leq x(t) \leq \beta(t)$ for every $t\in I$.
\begin{proof}
For each $\eta\in[\alpha(a),\beta(a)]$, consider the following Dirichlet problem
\begin{equation}\label{etaP}
    \left\{
    \begin{array}{lll}
   ((A \circ x)\cdot (\Phi \circ x'_{g_1}))'_{g_2}(t)\in G(t,x(t),x'_{g_1}(t)) \quad \mbox{for $g_2$-a.a. $t \in I=[a,b]$,}
   \\ x(a)=\eta, \quad x(b)=h(\eta), \\
   \alpha(t)\leq x(t) \leq \beta(t) \,\, \text{for every }\,\, t\in I.
   \end{array}
   \right.
\end{equation}
Since $h$ is increasing, from the definitions of lower and upper solutions we get 
$$
\alpha(b)=h(\alpha(a))\leq h(\eta) \leq h(\beta(a))=\beta(b),
$$
and since $\alpha(a)\leq\eta\leq\beta(a)$, $\alpha$ and $\beta$ are a well-ordered pair of lower and upper solutions for problem \eqref{etaP} as well. Thus, the requirements of Theorem \ref{teo:DP} are satisfied, so the set 
$$
S_{\eta}=\left\{ x\in\mathcal{BC}_{g_1}^1(I) \, : \, x \,\,\text{is a solution of \eqref{etaP}}\right\}
$$
is non-empty for every $\eta\in[\alpha(a),\beta(a)]$.

To prove the theorem, it suffices to find some $\eta^*\in[\alpha(a),\beta(a)]$ and $x\in S_{\eta^*}$ such that $$k(x(a),x(b),x'_{g_1}(a),x'_{g_1}(b))=0,$$ since in that case we would also have that $h(x(a))=h(\eta)=x(b)$ and $x$ would be a solution to problem \eqref{NLP} between $\alpha$ and $\beta$. To this end, assume by a contradiction that for every $\eta\in[\alpha(a),\beta(a)]$ and every $x\in S_{\eta}$ we have $k(x(a),x(b),x'_{g_1}(a),x'_{g_1}(b))\neq0$.

Observe now that if $\eta=\alpha(a)$ and $x\in S_{\alpha(a)}$, then $x(a)=\alpha(a)$ and $x(b)\geq\alpha(b)=h(\alpha(a))=h(\eta)=x(b)$, i.e., $x(b)=\alpha(b)$. Moreover, since $\alpha(t)\leq x(t)$ for every $t\in I$, we must have $x'_{g_1}(a)\geq\alpha'_{g_1}(a)$ and $x'_{g_1}(b)\leq\alpha'_{g_1}(b)$. Hence, from \eqref{condK} we deduce that
$$
k(x(a),x(b),x'_{g_1}(a),x'_{g_1}(b))\geq k(\alpha(a),\alpha(b),\alpha'_{g_1}(a),\alpha'_{g_1}(b))\geq 0, \quad \text{if}\,\, \eta=\alpha(a).
$$
By a completely analogous argument with $\eta=\beta(a)$, we get
$$
k(x(a),x(b),x'_{g_1}(a),x'_{g_1}(b))\leq k(\beta(a),\beta(b),\beta'_{g_1}(a),\beta'_{g_1}(b))\leq0, \quad \text{if}\,\, \eta=\beta(a).
$$
Let us define the set
$$
E=\left\{ \eta\in[\alpha(a),\beta(a)] \, : \, \text{there exists}\,\, x\in S_{\eta} \,\, \text{with}\,\, k(x(a),x(b),x'_{g_1}(a),x'_{g_1}(b))>0\right\}.
$$
We have that $\beta(a)\notin E$ and, by the contradiction assumption, $\alpha(a)\in E$. Now, let $\eta^*=\sup{E}$, take an increasing sequence $\{\eta_n\}_n\subset E$ converging to $\eta^*$ and let $\{x_n\}_n$ be a sequence such that $x_n\in S_{\eta_n}$ for each $n\in\mathbb{N}$. Since $h$ is continuous, $h(\eta_n)\to h(\eta^*)$, so applying Lemma \ref{lem:DC} we deduce the existence of a function $x^*\in S_{\eta^*}$ such that $x_n\to x^*$ in $\mathcal{BC}_{g_1}(I)$
and $(x_n)'_{g_1}(t) \to (x^*)'_{g_1}(t)$ for any $t \in I$ (passing to a subsequence if necessary). In particular, $(x_n)'_{g_1}(a)\to(x^*)'_{g_1}(a)$ and $(x_n)'_{g_1}(b)\to(x^*)'_{g_1}(b)$. Subsequently, since $k(x_n(a),x_n(b),(x_n)'_{g_1}(a),(x_n)'_{g_1}(b))>0$ for every $n\in \mathbb{N}$ and $k$ is continuous, from the assumption of contradiction we deduce that 
$$
k(x^*(a),x^*(b),(x^*)'_{g_1}(a),(x^*)'_{g_1}(b))>0,
$$
so $\eta^*\in E$ and $\eta^*<\beta(a)$.

Consider now the Dirichlet problem
\begin{equation}\label{etaP2}
    \left\{
    \begin{array}{lll}
   ((A \circ x)\cdot (\Phi \circ x'_{g_1}))'_{g_2}(t)\in G(t,x(t),x'_{g_1}(t)) \quad \mbox{for $g_2$-a.a. $t \in I=[a,b]$,}
   \\ x(a)=\eta, \quad x(b)=h(\eta), \\
   x^*(t)\leq x(t) \leq \beta(t) \,\, \text{for every }\,\, t\in I.
   \end{array}
   \right.
\end{equation}
with $\eta\in[\eta^*,\beta(a)]$. It is straightforward to observe that $x^*$ and $\beta$ form a well-ordered pair of lower and upper solutions for this problem and, since the hypotheses of Theorem \ref{teo:DP} are again satisfied, we deduce that the set 
$$
S^*_{\eta}=\left\{ x\in\mathcal{BC}^1_{g_1}(I) \, : \, x \,\,\text{is a solution of \eqref{etaP2}}\right\}
$$
is non-empty for every $\eta\in[\eta^*,\beta(a)]$. Now, let $\{\eta_k\}_k\subset(\eta^*,\beta(a)]$ be a decreasing sequence converging to $\eta^*$ and let $\{x_k\}_k$ be a sequence such that $x_k\in S^*_{\eta_k}$ for each $k\in\mathbb{N}$. Applying Lemma \ref{lem:DC} once again, we deduce the existence of a function $\widetilde{x}\in S^*_{\eta^*}$ such that, passing to a subsequence if necessary, $x_k\to\widetilde{x}$ in $\mathcal{BC}_{g_1}(I)$, $(x_n)'_{g_1}(a)\to(x^*)'_{g_1}(a)$ and $(x_n)'_{g_1}(b)\to(x^*)'_{g_1}(b)$. Since $\eta_k\notin E$, by the contradiction assumption we have that $k(x_k(a),x_k(b),(x_k)'_{g_1}(a),(x_k)'_{g_1}(b))<0$ for every $k\in \mathbb{N}$. Hence, using the same argument as above, we get 
\begin{equation}\label{kpos}
k(\widetilde{x}(a),\widetilde{x}(b),\widetilde{x}'_{g_1}(a),\widetilde{x}'_{g_1}(b))<0.
\end{equation}
Now, note that $x^*(t)\leq\widetilde{x}(t)$ for every $t\in I$, $x^*(a)=\widetilde{x}(a)$ and $x^*(b)=\widetilde{x}(b)$. Subsequently, $(x^*)'_{g_1}(a)\leq\widetilde{x}'_{g_1}(a)$ and $(x^*)'_{g_1}(b)\geq\widetilde{x}'_{g_1}(b)$, so by the monotonicity of $k$ we deduce that
$$
k(\widetilde{x}(a),\widetilde{x}(b),\widetilde{x}'_{g_1}(a),\widetilde{x}'_{g_1}(b))\geq k(x^*(a),x^*(b),(x^*)'_{g_1}(a),(x^*)'_{g_1}(b))>0,
$$
a contradiction with \eqref{kpos}. Thus, there must exist some $\eta^*\in[\alpha(a),\beta(a)]$ and $x\in S_{\eta^*}$ with $$
k(x(a),x(b),x'_{g_1}(a),x'_{g_1}(b))=0,$$
which means that $x$ is a solution for \eqref{NLP}.
\end{proof}
\end{theorem}

It is worth noting that Theorem \ref{teo:NLP} guarantees the existence of a solution for the periodic problem
\begin{equation}\label{PP}
    ((A \circ x)\cdot (\Phi \circ x'_{g_1}))'_{g_2}(t)\in G(t,x(t),x'_{g_1}(t)) \quad \mbox{for $g_2$-a.a. $t \in I=[a,b]$,} \quad x(a)=x(b),\quad x'_{g_1}(a)=x'_{g_1}(b),
\end{equation}
Indeed, by setting $h(x)=x$ and $k(x,y,z,w)=z-w$ in \eqref{NLP} we obtain the following result, taking into account that the $g_1$-slope in problem \eqref{PP} is null.

\begin{corollary}
     In the conditions of Theorem \ref{teo:DP} with $\alpha$ and $\beta$   a well-ordered pair of lower and upper solutions for problem \eqref{PP}, assume that condition \eqref{condK} holds and assume also that 
     $$
     M|\Phi(0)| <cm-2\|\psi\|_{L^{1}_{g_2}(I)}.
     $$
 \noindent
 Then, problem \eqref{PP} has at least one solution $x(t)$ such that $\alpha(t)\leq x(t) \leq \beta(t)$ for every $t\in I$.
\end{corollary}
We will now move onto the following problem with a different type of nonlinear boundary conditions
\begin{equation}\label{NLP2}
    \left\{
    \begin{array}{lll}
  ((A \circ x)\cdot (\Phi \circ x'_{g_1}))'_{g_2}(t)\in G(t,x(t),x'_{g_1}(t)) \quad \mbox{for $g_2$-a.a. $t \in I=[a,b]$,}
   \\ p(x(a),x'_{g_1}(a))=0,\\
   q(x(b),x'_{g_1}(b))=0,
   \end{array}
   \right.
\end{equation}
where $p,q:\mathbb{R}^2\rightarrow\mathbb{R}$ are continuous functions such that for each fixed $x$ we have
\begin{equation}\label{condPQ}
y\to p(x,y)\,\,\, \text {is decreasing}\quad \text{and} \quad y\to q(x,y) \, \, \,\text{is increasing}.  
\end{equation}

\noindent
As expected, lower and upper solutions for the problem \eqref{NLP2} are defined in the following way.
\begin{definition}
    A lower solution for the problem \eqref{NLP2} is a function $\alpha\in\mathcal{BC}^1_{g_1}(I)$ such that 
    $$p(\alpha(a),\alpha'_{g_1}(a))\leq 0,\quad q(\alpha(b),\alpha'_{g_1}(b))\leq 0,
    $$
    and $(A \circ \alpha)\cdot (\Phi \circ \alpha'_{g_1}) \in\mathcal{AC}_{g_2}(I)$ and
   \begin{equation*}
   G(t,\alpha(t),\alpha'_{g_1}(t))\cap\left(-\infty,  ((A \circ \alpha)\cdot (\Phi \circ \alpha'_{g_1}))'_{g_2}(t)\right]\neq\emptyset \quad \text{for}\,\, g_2\text{-almost}\,\, t \in I.
   \end{equation*}

   Analogously, an upper solution for the problem \eqref{NLP2} is a function $\beta\in\mathcal{BC}^1_{g_1}(I)$ such that 
    $$
    p(\beta(a),\beta'_{g_1}(a))\geq 0,\quad q(\beta(b),\beta'_{g_1}(b))\geq 0,
    $$
    and $(A \circ \beta)\cdot (\Phi \circ \beta'_{g_1}) \in\mathcal{AC}_{g_2}(I)$ and 
   \begin{equation*}
   G(t,\beta(t),\beta'_{g_1}(t))\cap\left[  ((A \circ \beta)\cdot (\Phi \circ \beta'_{g_1})'_{g_2}(t),+\infty\right)\neq\emptyset \quad \text{for}\,\, g_2\text{-almost}\,\, t \in I.
   \end{equation*}
\end{definition}

Relying again on the existence theorem for the Dirichlet problem, we get to the following result.
\begin{theorem}\label{teo:NLP2}
In the conditions of Theorem \ref{teo:DP} with $\alpha$ and $\beta$   a well-ordered pair of lower and upper solutions for problem \eqref{NLP2}, assume that condition \eqref{condPQ} holds and assume also that  
 \begin{equation}\label{ms2}
     M\Phi\left(\dfrac{\beta(b)-\alpha(a)}{g_1(b)-g_1(a)}\right)<cm-2\|\psi\|_{L_{g_2}}, \quad  M\Phi\left(\dfrac{\beta(a)-\alpha(b)}{g_1(b)-g_1(a)}\right)>2\|\psi\|_{L_{g_2}}-cm.
 \end{equation}
 Then, problem \eqref{NLP2} has at least one solution $x(t)$ such that $\alpha(t)\leq x(t) \leq \beta(t)$ for every $t\in I$.
\begin{proof}
For each $(\eta,\mu)\in R=[\alpha(a),\beta(a)]\times[\alpha(b),\beta(b)]$, consider the Dirichlet problem
\begin{equation}\label{muP}
    \left\{
    \begin{array}{lll}
   ((A \circ x)\cdot (\Phi \circ x'_{g_1}))'_{g_2}\in G(t,x(t),x'_{g_1}(t)) \quad \mbox{for $g_2$-a.a. $t \in I=[a,b]$,}
   \\ x(a)=\eta, \quad x(b)=\mu, \\
   \alpha(t)\leq x(t) \leq \beta(t) \,\, \text{for every }\,\, t\in I.
   \end{array}
   \right.
\end{equation}
    For every $(\eta,\mu)\in R$, we have that
    $$
    \dfrac{\alpha(b)-\beta(a)}{g_1(b)-g_1(a)}\leq\frac{\mu-\eta}{g_1(b)-g_1(a)}\leq\dfrac{\beta(b)-\alpha(a)}{g_1(b)-g_1(a)},
    $$
    so, since $\Phi$ is increasing, assumption \eqref{ms2} ensures that condition \eqref{ms} is satisfied for every $s^*$ of the form $(\mu-\eta)/(g_1(b)-g_1(a))$ with $(\eta,\mu)\in R$. Hence, we can apply Theorem \ref{teo:DP} to deduce existence of solutions for problem \eqref{muP} between $\alpha$ and $\beta$. Moreover, from the definition of lower and upper solutions and the monotonicity of $p$, for every $x(t)$ solution of \eqref{muP} such that $\alpha(t)\leq x(t)\leq\beta(t)$ for every $t\in I$, we have
    \begin{equation}\label{palpha}
    \begin{split}
    p(x(a),x'_{g_1}(a))\leq p(\eta,\alpha'_{g_1}(a))\leq0, \quad \text{if}\,\,\eta=\alpha(a),\\
    p(x(a),x'_{g_1}(a))\geq p(\eta,\beta'_{g_1}(a))\leq0, \quad \text{if}\,\,\eta=\beta(a),
    \end{split}
    \end{equation}
    since $x'_{g_1}(a)\geq\alpha'_{g_1}(a)$ if $x(a)=\alpha(a)$, and $x'_{g_1}(a)\leq\beta'_{g_1}(a)$ if $x(a)=\beta(a)$. Following the same reasoning as in the proof of Theorem \ref{teo:NLP}, specifically regarding the existence of a solution $x\in\mathcal{BC}^1_{g_1}(I)$ with an initial value $\eta\in[\alpha(a),\beta(a)]$ satisfying $k(x(a),x(b),x'_{g_1}(a),x'_{g_1}(b))=0$, one can show that for every $\mu\in[\alpha(b),\beta(b)]$, there exists some $\eta\in[\alpha(a),\beta(a)]$ such that the following problem
    \begin{equation}\label{muP2}
    \left\{
    \begin{array}{lll}
    ((A \circ x)\cdot (\Phi \circ x'_{g_1}))'_{g_2}\in G(t,x(t),x'_{g_1}(t)) \quad \mbox{for $g_2$-a.a. $t \in I=[a,b]$,}\\ 
   x(a)=\eta, \quad x(b)=\mu, \\
   p(x(a),x'_{g_1}(a))=0, \\
   \alpha(t)\leq x(t) \leq \beta(t) \,\, \text{for every }\,\, t\in I,
   \end{array}
   \right.
    \end{equation}
    admits solutions. Subsequently, if for each $\mu\in[\alpha(b),\beta(b)]$ we define the sets
    $$
    E_{\mu}=\left\{ \eta\in[\alpha(a),\beta(a)] \, : \, \text{\eqref{muP2} has solutions}\right\},
    $$
    $$
    S_{\mu,\eta}=\{x\in\mathcal{AC}_{g_1}(I) \, : \, x \,\,\text{is a solution of \eqref{muP2}}\},
    $$
    we have that $E_{\mu}\neq\emptyset$ for every $\mu\in[\alpha(b),\beta(b)]$ and $S_{\mu,\eta}\neq\emptyset$ for every $(\mu,\eta)\in R$. 

    To prove the result, it suffices to see that there exist some $\mu^*\in[\alpha(b),\beta(b)]$, $\eta^*\in E_{\mu^*}$ and $x^*\in S_{\mu^*,\eta^*}$ such that $q(x^*(b),(x^*)'_{g_1}(b))=0$, since in that case $x^*$ would be a solution for problem \eqref{NLP2} between $\alpha$ and $\beta$.

    \noindent
    To this end, assume by a contradiction that for every $\mu\in[\alpha(b),\beta(b)]$, every $\eta\in E_{\mu}$ and every $x\in S_{\mu,\eta}$ we have that $q(x(b),x'_{g_1}(b))\neq0$. Observe that, by the monotonicity of $q$, every $x\in S_{\mu,\eta}$ verifies 
    \begin{equation}\label{qmu}
        \begin{split}
            q(x(b),x'_{g_1}(b))\leq q(\mu,\alpha'_{g_1}(b))\leq 0, \quad \text{if}\,\, \mu=\alpha(b), \\
            q(x(b),x'_{g_1}(b))\geq q(\mu,\beta'_{g_1}(b))\geq 0, \quad \text{if}\,\, \mu=\alpha(b),
        \end{split}
    \end{equation}
    since $x'_{g_1}(b)\leq\alpha'_{g_1}(b)$ if $x(b)=\alpha(b)$, and $x'_{g_1}(b)\geq\beta'_{g_1}(b)$ if $x(b)=\beta(b)$. 

    \noindent
    Let
    $$
    N=\{\mu\in[\alpha(b),\beta(b)]\,:\,q(x(b),x'_{g_1}(b))<0\,\, \text{for some}\,\, x\in S_{\mu,\eta}, \eta\in E_{\mu}\}.
    $$
    By \eqref{qmu} and the contradiction assumption, we deduce that $\alpha(b)\in N$ and $\beta(b)\notin N$. Let $\mu^*=\sup{N}$, $\{\mu_n\}_n\subset E$ and increasing sequence converging to $\mu^*$, $\{\eta_n\}_n\subset [\alpha(a),\beta(a)]$ such that $\eta_n\in E_{\mu_n}$ for every $n\in\mathbb{N}$ and $\{x_n\}_n\subset\mathcal{AC}_{g_1}(I)$ such that $x_n\in S_{\mu_n,\eta_n}$ and 
    \begin{equation}\label{qneg}
        q(x_n(b),(x_n)'_{g_1}(b))<0,
    \end{equation}
    for every $n\in\mathbb{N}$. Since $[\alpha(a),\beta(a)]$ is compact, $\{\eta_n\}_n$ has a convergent subsequence, relabeled again as $\{\eta_n\}_n$, such that $\eta_n\to\eta^*\in[\alpha(a),\beta(a)]$. Now, applying Lemma \ref{lem:DC}, we deduce the existence of a function $x^*\in S_{\mu^*,\eta^*}$ such that, passing to a subsequence if necessary, $x_n\to x^*$ in $\mathcal{BC}_{g_1}(I)$ and $q(x_n(b),(x_n)'_{g_1}(b))\to q(x^*(b),(x^*)'_{g_1}(b))$. By \eqref{qneg} and the contradiction assumption, we get that $q(x^*(b),(x^*)'_{g_1}(b))<0$. Hence, $\mu^*\in E$ and $\mu^*<\beta(b)$.

    Consider now the Dirichlet problem
     \begin{equation}\label{muP3}
    \left\{
    \begin{array}{lll}
   ((A \circ x)\cdot (\Phi \circ x'_{g_1}))'_{g_2}\in G(t,x(t),x'_{g_1}(t)) \quad \mbox{for $g_2$-a.a. $t \in I=[a,b]$,}\\ 
   x(a)=\eta, \quad x(b)=\mu, \\
   p(x(a),x'_{g_1}(a))=0, \\
   x^*(t)\leq x(t) \leq \beta(t) \,\, \text{for every }\,\, t\in I.
   \end{array}
   \right.
    \end{equation}
     Note that $x^*$ and $\beta$ form a well-ordered pair of lower and upper solutions for \eqref{muP3}, so by the same reasoning as before (replacing $\alpha$ by $x^*$), we get that the sets 
     $$
    E^*_{\mu}=\left\{ \eta\in[x^*(a),\beta(a)] \, : \, \text{\eqref{muP3} has solutions}\right\},
    $$
    $$
    S^*_{\mu,\eta}=\{x\in\mathcal{AC}_{g_1}(I) \, : \, x \,\,\text{is a solution of \eqref{muP3}}\},
    $$
    are non-empty for every $\mu\in[x^*(b),\beta(b)]$ and every $\eta\in E^*_\mu$, respectively.

    \noindent
    Since $\mu^*<\beta(b)$, we can take a decreasing sequence $\{\mu_m\}_m\subset(\mu^*,\beta(b))$ such that $\mu_m\to\mu^*$ and, subsequently, there exist $\{\eta_m\}_m\subset[\eta^*,\beta(a)]$ and $\{x_m\}_m\subset\mathcal{BC}^1_{g_1}(I)$ such that $\eta_m\in E^*_{\mu_m}$, $x_m\in S^*_{\mu_m,\eta_m}$ and
    \begin{equation}\label{qpos}
        q(x_m(b),(x_m)'_{g_1}(b))>0
    \end{equation}
    for every $m\in\mathbb{N}$, since $\mu_m\notin N$. Applying Lemma \ref{lem:DC} once again and passing to subsequences if necessary, we deduce the existence if some $\widetilde{\eta}\in[\eta^*,\beta(a)]$ such that $\eta_m\to\widetilde{\eta}$ and some $\widetilde{x}\in S^*_{\mu^*,\tilde{\eta}}$ such that $x_m\to\widetilde{x}$ in $\mathcal{BC}_{g_1}(I)$ and $(x_m)'_{g_1}\to\widetilde{x}'_{g_1}$ pointwise. Since $q(x_m(b),(x_m)'_{g_1}(b))\to q(\widetilde{x}(b),\widetilde{x}'_{g_1}(b))$, from \eqref{qpos} and the contradiction assumption we get 
    \begin{equation}\label{qpos2}
    q(\widetilde{x}(b),\widetilde{x}'_{g_1}(b))>0.
    \end{equation}
    Now, since $x^*(t)\leq\widetilde{x}(t)$ for every $t\in I$ and $x^*(b)=\mu^*=\widetilde{x}(b)$, necessarily $(x^*)'_{g_1}(b)\geq \widetilde{x}'_{g_1}(b)$, and from the monotonicity of $q$ we deduce that 
    $$
    q(x^*(b),\widetilde{x}'_{g_1}(b))\leq q(\widetilde{x}(b),(x^*)'_{g_1}(b))<0,
    $$
    a contradiction with \eqref{qpos2}. Thus, there must exist some $\mu^*\in[\alpha(b),\beta(b)]$, $\eta^*\in E_{\mu^*}$ and $x\in S_{\mu^*,\eta^*}$ with $q(x^*(b),(x^*)'_{g_1}(b))=0$, which means that $x^*$ is a solution for \eqref{NLP2}.
\end{proof}
\end{theorem}

\begin{remark}
    As we mentioned at the beginning of the Section, the nonlinear formulation of the boundary conditions in \eqref{NLP2} covers some particular cases of interest.
    
    \noindent
    On the one hand, taking $p(x,y)=c_ax-d_ay-v_a$ and $q(x,y)=c_bx+d_by-v_b$, with $d_a, \, d_b \ge 0$ and $(|c_a|+d_a)(|c_b|+d_b)>0$,
    on Theorem \ref{teo:NLP2} we deduce the existence of a solution of the problem with Robin boundary conditions
    $$
    \left\{
    \begin{array}{l}
     ((A \circ x)\cdot (\Phi \circ x'_{g_1}))'_{g_2}\in G(t,x(t),x'_{g_1}(t)) \quad \mbox{for $g_2$-a.a. $t \in I=[a,b]$,}\\
     \\
     c_a x(a)-d_a x'(a)=v_a, \,\, c_b x(b)+d_b x'(b)=v_b.
     \end{array}
     \right.
    $$
    On the other hand, taking $p(x,y)=v_a-y$ and $q(x,y)=y-v_b$, we obtain the existence of a solution for the analogous problem with Neumann boundary conditions $x'_{g_1}(a)=v_a$ and $x'_{g_1}(b)=v_b$.
    
\end{remark}


\end{document}